\documentclass[leqno,11pt,a4paper,twoside]{article}%
\usepackage{amssymb}
\usepackage{amsfonts}
\usepackage{amsmath}
\usepackage{graphicx}%
\providecommand{\U}[1]{\protect\rule{.1in}{.1in}}
\newtheorem{theorem}{Theorem}

\newtheorem{criterion}[theorem]{Criterion}
\newtheorem{definition}[theorem]{Definition}

\newtheorem{proposition}[theorem]{Proposition}
\newtheorem{remark}[theorem]{Remark}

\newenvironment{proof}[1][Proof]{\noindent\textbf{#1.} }{\ \rule{0.5em}{0.5em}}
\begin{document}

\title{Viability, Invariance and Reachability for Controlled Piecewise Deterministic
Markov Processes Associated to Gene Networks}
\author{Goreac, D.\thanks{Universit\'{e} Paris-Est, Laboratoire d'Analyse et
Math\'{e}matiques Appliqu\'{e}es, UMR 8050, Boulevard Descartes, Cit\'{e}
Descartes, 77450, Champs-sur-Marne}}
\maketitle

\begin{abstract}
We aim at characterizing viability, invariance and some reachability
properties of controlled piecewise deterministic Markov processes (PDMPs).
Using analytical methods from the theory of viscosity solutions, we establish
criteria for viability and invariance in terms of the first order normal cone.
We also investigate reachability of arbitrary open sets. The method is based
on viscosity techniques and duality for some associated linearized problem.
The theoretical results are applied to general On/Off systems, Cook's model
for haploinssuficiency, and a stochastic model for bacteriophage $\lambda$.

\end{abstract}



\section{Introduction}

Markov processes have been intensively used to describe variability features
of various cellular processes. To our best knowledge, Markovian tools have
first been employed in connection to molecular biology in \cite{De}. The
natural idea was to associate to each reaction network a pure jump model. Due
to the large number of molecular species involved in the reactions, direct
simulation of these models turns out to be very slow. To increase proficiency,
hybrid models are adopted in \cite{CDR}. They distinguish the discrete
components from the "continuous" ones. Using partial Kramers-Moyal expansion,
the authors of \cite{CDR} replace the initial pure jump process with an
appropriate piecewise deterministic Markov one.

One may reduce the complexity of PDMPs by restricting the model to some
invariant set containing the initial data, whenever this is known. Compact
invariant sets are also needed for efficiently implementing algorithms.
Another important issue that can be approached using invariance are the stable
points. In particular, a fixed point for which one finds arbitrarily small
surrounding invariant sets is stable in the sense of Lyapunov.

We begin by characterizing $\varepsilon$-viability of controlled PDMPs via
some associated control problem. A closed set of constraints $K$ is said to be
viable (or $\varepsilon$-viable) with respect to some dynamic control system
if, starting from $K,$ one is able to find suitable controls keeping the
trajectory in $K$ (or, at least in some arbitrarily small neighborhood of the
set of constraints). Viability properties have been extensively studied in
both deterministic and stochastic settings (for Brownian diffusions), starting
from the pioneer work of Nagumo. The methods used to describe this property
for deterministic or diffusion processes rely either on the Bouligand-Severi
contingent cone (cf. \cite{A}, \cite{ADP}, \cite{GT}) or on viscosity
solutions (\cite{AF}, \cite{BC}, \cite{BG}, \cite{BJ}, \cite{BPQR}). Using
analytical tools from viscosity theory, we provide a geometrical
characterization of $\varepsilon$-viability and invariance of some set of
constraints $K$ with respect to the controlled piecewise deterministic Markov
process. As for the Brownian diffusion case (cf. \cite{BJ}), the criterion
involves the normal cone to the set of constraints and is completely
deterministic. Similar arguments allow one to characterize the invariance of
the set of constraints. We emphasize that these geometrical conditions can be
rather easily checked for PDMPs associated to gene networks. In order to
illustrate these theoretical assertions, two examples are considered. For
general On/Off models, we show how the invariance criterion can be used in
order to reduce the state space to a compact set. We also characterize points
that can be chosen as candidates for stability (in the sense that one finds
arbitrarily small surrounding regions that are invariant). Another biological
example is a model for bacteriophage $\lambda$ (described in \cite{HPDC}).
Although it is more complex, one can still use the invariance criterion to
characterize candidates for stability. It turns out that only one such point
exists in the absence of impulsive exterior control factors.

The second aim of the paper is to characterize the reachability property of
arbitrary open sets with respect to the controlled piecewise deterministic
Markov process. The criterion is obtained using viscosity methods. Recently,
the paper \cite{BGQ} has provided a linear programming formulation for
discounted control problems in the framework of SDEs driven by standard
Brownian motion. The reachability problem can be connected to the value
function of some appropriate piecewise deterministic control system. Using the
idea in \cite{BGQ}, we give a criterion involving the dual formulation of the
linearized version of the initial problem. To illustrate this result, we
consider Cook's model for haploinsufficiency introduced in \cite{CGT}. Our
criterion allows one to prove that, starting from any arbitrary point, one
reaches any arbitrarily given open region, with positive probability.

The paper is organized as follows: In Subsection \ref{constr+ass} we briefly
recall the construction of controlled PDMPs and state the main (standard)
assumptions. Section \ref{Section1} is devoted to the study of viability
property (Subsection \ref{Subsection1.1}) and invariance (Subsection
\ref{Subsection1.2}) with respect to the PDMP. The criteria involve the normal
cone to the set of constraints and the characteristics of the process. Section
\ref{Section2} deals with the reachability property. We use a Krylov-type
argument to provide some dual formulation of the associated control problem.
In Subsection \ref{Subsection3.1} we recall some rudiments on the PDMPs
associated to a system of chemical reactions. We consider two biological
examples: the On/Off model (Subsection \ref{Subsection3.2}) and the
bacteriophage $\lambda$ (Subsection \ref{Subsection3.3.}). We first study the
compact invariant sets for the On/Off model. For a particular case (the
so-called Cook model for haploinsufficiency), we prove that every open set can
be reached with positive probability, starting from any initial point. In the
case of the bacteriophage $\lambda$ (described in \cite{HPDC}), our invariance
criterion allow to identify the stable point of the system. The Appendix
provides the comparison principle and some stability results for viscosity solutions.

\subsection{Construction of controlled PDMPs and main
assumptions\label{constr+ass}}

We let $U$ be a compact metric space (the control space) and $%
\mathbb{R}
^{N}$ be the state space, for some $N\geq1.$

Piecewise deterministic control processes have been introduced by Davis
\cite{D}. Such processes are given by their local characteristics: a vector
field $f:%
\mathbb{R}
^{N}\times U\rightarrow%
\mathbb{R}
^{N}$ that determines the motion between two consecutive jumps, a jump rate
$\lambda:%
\mathbb{R}
^{N}\times U\rightarrow%
\mathbb{R}
_{+}$ and a transition measure $Q:%
\mathbb{R}
^{N}\times U\times\mathcal{B}\left(
\mathbb{R}
^{N}\right)  \rightarrow\mathcal{P}\left(
\mathbb{R}
^{N}\right)  .$ Here $\mathcal{B}\left(
\mathbb{R}
^{N}\right)  $ is the family of Borel sets and $\mathcal{P}\left(
\mathbb{R}
^{N}\right)  $ stands for the family of probability measures on $%
\mathbb{R}
^{N}.$ For every $A\in\mathcal{B}\left(
\mathbb{R}
^{N}\right)  ,$ the function $\left(  u,x\right)  \mapsto Q\left(
x,u,A\right)  $ should be measurable and, for every $\left(  x,u\right)  \in%
\mathbb{R}
^{N}\times U$, $Q\left(  x,u,\left\{  x\right\}  \right)  =0.$

We summarize the construction of controlled piecewise deterministic Markov
processes (PDMP). Whenever $u\in\mathbb{L}^{0}\left(
\mathbb{R}
^{N}\times%
\mathbb{R}
_{+};U\right)  $ ($u$ is a Borel measurable function) and $\left(  t_{0}%
,x_{0}\right)  \in%
\mathbb{R}
_{+}\times%
\mathbb{R}
^{N},$ we consider the ordinary differential equation%
\[
\left\{
\begin{array}
[c]{l}%
d\Phi_{t}^{t_{0},x_{0},u}=f\left(  \Phi_{t}^{t_{0},x_{0},u},u\left(
x_{0},t-t_{0}\right)  \right)  dt,\text{ }t\geq t_{0},\\
\Phi_{t_{0}}^{t_{0},x_{0},u}=x_{0.}%
\end{array}
\right.
\]
We choose the first jump time $T_{1}$ such that the jump rate be
$\lambda\left(  \Phi_{t}^{0,x_{0},u},u\left(  x_{0},t\right)  \right)  $%
\[
\mathbb{P}\left(  T_{1}\geq t\right)  =\exp\left(  -\int_{0}^{t}\lambda\left(
\Phi_{s}^{0,x_{0},u},u\left(  x_{0},s\right)  \right)  ds\right)  .
\]
The controlled piecewise deterministic Markov processes (PDMP) is defined by
\[
X_{t}^{x_{0},u}=\Phi_{t}^{0,x_{0},u},\text{ if }t\in\left[  0,T_{1}\right)  .
\]
The post-jump location $Y_{1}$ has $Q\left(  \Phi_{\tau}^{0,x_{0},u},u\left(
x_{0},\tau\right)  ,\cdot\right)  $ as conditional distribution given
$T_{1}=\tau.$ Starting from $Y_{1}$ at time $T_{1}$, we select the inter-jump
time $T_{2}-T_{1}$ such that
\[
\mathbb{P}\left(  T_{2}-T_{1}\geq t\text{ }/\text{ }T_{1},Y_{1}\right)
=\exp\left(  -\int_{T_{1}}^{T_{1}+t}\lambda\left(  \Phi_{s}^{T_{1},Y_{1}%
,u},u\left(  Y_{1},s-T_{1}\right)  \right)  ds\right)  .
\]
We set
\[
X_{t}^{x_{0},u}=\Phi_{t}^{T_{1},Y_{1},u},\text{ if }t\in\left[  T_{1}%
,T_{2}\right)  .
\]
The post-jump location $Y_{2}$ satisfies%
\[
\mathbb{P}\left(  Y_{2}\in A\text{ }/\text{ }T_{1},Y_{1}\right)  =Q\left(
\Phi_{T_{2}}^{T_{1},Y_{1},u},u\left(  Y_{1},T_{2}-T_{1}\right)  ,A\right)  ,
\]
for all Borel set $A\subset%
\mathbb{R}
^{N}.$ And so on.

Throughout the paper, unless stated otherwise, we assume the following:

(A1) The function $f:%
\mathbb{R}
^{N}\times U\longrightarrow%
\mathbb{R}
^{N}$ is uniformly continuous on $%
\mathbb{R}
^{N}\times U$ and there exists a positive real constant $C>0$ such that
\begin{equation}
\left\vert f\left(  x,u\right)  -f\left(  y,u\right)  \right\vert \leq
C\left\vert x-y\right\vert ,\text{ and }\left\vert f\left(  x,u\right)
\right\vert \leq C, \tag{A1}\label{A1}%
\end{equation}
for all $x,y\in%
\mathbb{R}
^{N}$ and all $u\in U.$

(A2) The function $\lambda:%
\mathbb{R}
^{N}\times U\longrightarrow%
\mathbb{R}
_{+}$ is uniformly continuous on $%
\mathbb{R}
^{N}\times U$ and there exists a positive real constant $C>0$ such that
\begin{equation}
\left\vert \lambda\left(  x,u\right)  -\lambda\left(  y,u\right)  \right\vert
\leq C\left\vert x-y\right\vert ,\text{ and }\lambda\left(  x,u\right)  \leq
C, \tag{A2}\label{A2}%
\end{equation}
for all $x,y\in%
\mathbb{R}
^{N}$ and all $u\in U.$

(A3) For each bounded uniformly continuous function $h\in BUC\left(
\mathbb{R}
^{N}\right)  ,$ there exists a continuous function $\eta_{h}:%
\mathbb{R}
\longrightarrow%
\mathbb{R}
$ such that $\eta_{h}\left(  0\right)  =0$ and
\begin{equation}
\sup_{u\in U}\left\vert \int_{%
\mathbb{R}
^{N}}h\left(  z\right)  Q\left(  x,u,dz\right)  -\int_{%
\mathbb{R}
^{N}}h\left(  z\right)  Q\left(  y,u,dz\right)  \right\vert \leq\eta
_{h}\left(  \left\vert x-y\right\vert \right)  . \tag{A3}\label{A3}%
\end{equation}

(A4) For every $x\in%
\mathbb{R}
^{N}$ and every decreasing sequence $\left(  \Gamma_{n}\right)  _{n\geq0}$ of
subsets of $%
\mathbb{R}
^{N},$%
\begin{equation}
\inf_{n\geq0}\sup_{u\in U}Q\left(  x,u,\Gamma_{n}\right)  =\sup_{u\in
U}Q\left(  x,u,\underset{n}{\cap}\Gamma_{n}\right)  . \tag{A4}\label{A4}%
\end{equation}

\begin{remark}
\label{rem1}We have kept A3 as it appears in Soner \cite{S}. However, this
assumption may be somewhat weakened by imposing

(A3') For each bounded uniformly continuous function $h\in BUC\left(
\mathbb{R}
^{N}\right)  ,$ there exists a continuous function $\eta_{h}:%
\mathbb{R}
\longrightarrow%
\mathbb{R}
$ such that $\eta_{h}\left(  0\right)  =0$ and
\[
\sup_{u\in U}\left\vert \lambda\left(  x,u\right)  \int_{%
\mathbb{R}
^{N}}h\left(  z\right)  Q\left(  x,u,dz\right)  -\lambda\left(  y,u\right)
\int_{%
\mathbb{R}
^{N}}h\left(  z\right)  Q\left(  y,u,dz\right)  \right\vert \leq\eta
_{h}\left(  \left\vert x-y\right\vert \right)  .
\]
It is obvious that whenever one assumes A3 and $\lambda\left(  \cdot\right)  $
is bounded, the assumption A3' \ holds true. Moreover, all the proofs in this
paper can be obtained (with minor changes) when A3' replaces A3.
\end{remark}

\section{A geometric condition for viability and invariance\label{Section1}}

\subsection{Conditions for viability\label{Subsection1.1}}

This subsection aims at characterizing the viability property of a nonempty,
closed set $K\subset%
\mathbb{R}
^{N}.$ In analogy to the deterministic framework, this property is proved to
be connected to some geometric condition involving the normal cone to $K.$ The
proof of the characterization relies on the viscosity solution concept. We
begin the subsection by recalling the notions of viability (respectively
$\varepsilon$-viability) and normal cone.

\begin{definition}
1. A nonempty, closed set $K\subset%
\mathbb{R}
^{N}$ is said to be viable with respect to the controlled piecewise
deterministic Markov process $X$ if, for every initial point $x\in K$, there
exists an admissible control process $u\in\mathbb{L}^{0}\left(
\mathbb{R}
^{N}\times%
\mathbb{R}
_{+};U\right)  $ such that $X_{t}^{x,u}\in K$, $\mathbb{P}$-almost surely, for
all $t\geq0.$

2. A nonempty, closed set $K\subset%
\mathbb{R}
^{N}$ is said to be $\varepsilon$-viable with respect to the controlled
piecewise deterministic process $X$ if, for every initial point $x\in K$ and
every $\varepsilon>0$, there exists an admissible control process
$u^{\varepsilon}\in\mathbb{L}^{0}\left(
\mathbb{R}
^{N}\times%
\mathbb{R}
_{+};U\right)  $ such that
\[
\mathbb{E}\left[  \int_{0}^{\infty}e^{-t}\left(  d_{K}\left(  X_{t}%
^{x,u^{\varepsilon}}\right)  \wedge1\right)  dt\right]  \leq\varepsilon.
\]
Here, $d_{K}$ stands for the distance function to the closed set $K$.
\end{definition}

\begin{definition}
Let $K\subset%
\mathbb{R}
^{N}$ be a closed subset and let $x$ be a point of $K.$ The normal cone to $K$
at $x,$ denoted by $N_{K}\left(  x\right)  $, is defined as%
\[
N_{K}\left(  x\right)  =\left\{  p\in%
\mathbb{R}
^{N}:\forall\varepsilon>0,\exists\eta>0\text{ such that }\forall y\in
K\cap\overline{B}\left(  x,\eta\right)  ,\text{ }\left\langle
p,y-x\right\rangle \leq\varepsilon\left\vert y-x\right\vert \right\}  .
\]
We recall that $\overline{B}\left(  x,\eta\right)  =\left\{  y\in%
\mathbb{R}
^{N}:\left\vert y-x\right\vert \leq\eta\right\}  .$
\end{definition}

The definition of the $\varepsilon$-viability property of a nonempty, closed
set $K\subset%
\mathbb{R}
^{N}$ can, alternatively, be given with respect to the value function%
\begin{equation}
v(x)=\inf_{u\in\mathbb{L}^{0}\left(
\mathbb{R}
^{N}\times%
\mathbb{R}
_{+};U\right)  }\mathbb{E}\left[  \int_{0}^{\infty}e^{-t}\left(  d_{K}\left(
X_{t}^{x,u}\right)  \wedge1\right)  dt\right]  , \label{v}%
\end{equation}
for all $x\in%
\mathbb{R}
^{N}.$ Indeed, with this notation, the set $K$ is $\varepsilon$-viable if and
only if the restriction of $v$ to $K$ is zero. We consider the associated
Hamilton-Jacobi integro-differential equation%
\begin{equation}
v\left(  x\right)  -d_{K}\left(  x\right)  \wedge1+H\left(  x,\nabla v\left(
x\right)  ,v\right)  =0, \label{HJB}%
\end{equation}
for all $x\in%
\mathbb{R}
^{N}$, where the Hamiltonian is given by
\begin{equation}
H\left(  x,p,\psi\right)  =\sup_{u\in U}\left\{  -\left\langle f\left(
x,u\right)  ,p\right\rangle -\lambda\left(  x,u\right)  \int_{%
\mathbb{R}
^{n}}\left(  \psi\left(  z\right)  -\psi\left(  x\right)  \right)  Q\left(
x,u,dz\right)  \right\}  . \label{H}%
\end{equation}
Under the assumptions (A1)-(A3), the function $v$ is known to satisfy (cf.
\cite{S}, Theorem 1.1), in the viscosity sense, Equation (\ref{HJB}). We are
going to need a slightly more general definition for the viscosity subsolution
(respectively supersolution) then the one used in \cite{S}.

\begin{definition}
A bounded, upper (lower) semicontinuous function $v$ is a viscosity
subsolution (supersolution) of (\ref{HJB}) if, for any test-function
$\varphi\in C_{b}^{1}\left(  N_{x}\right)  $, on some neighborhood $N_{x}$ of
$x\in%
\mathbb{R}
^{N}$, whenever $x$ is a maximum (minimum) point of $v-\varphi$,
\[
v\left(  x\right)  -d_{K}\left(  x\right)  \wedge1+H\left(  x,\nabla
\varphi\left(  x\right)  ,v\right)  \leq0\text{ (}\geq0\text{)}.
\]

A bounded, continuous function $v$ is a viscosity solution of (\ref{HJB}) if
it is both subsolution and supersolution.
\end{definition}

At this point, we introduce a technical assumption on the transition measure
$Q$ which provides a comparison principle. It states that the probability for
the post jump position to be arbitrarily far away from the pre jump one is
uniformly small. We emphasize that this assumption is made in order to give a
simple proof for the comparison principle. However, it is not essential; one
can, as an alternative, strengthen A3 as in \cite{AT} Section 3. Moreover, the
main results of the paper hold true independently of this assumption, whenever
a comparison principle for semicontinuous functions holds true.

(A5) We assume that
\begin{equation}
\inf_{n\geq1}\sup_{x\in%
\mathbb{R}
^{N},u\in U}Q\left(  x,u,%
\mathbb{R}
^{N}\smallsetminus\overline{B}\left(  x,n\right)  \right)  =0.\tag{A5}%
\label{A5}%
\end{equation}

\begin{remark}
Assumption (A5) is obviously satisfied if $Q$ does not depend on $x$ and $u$.
Moreover, all the piecewise deterministic processes associated to chemical
reactions (see Subsection \ref{Subsection3.1}) satisfy (A5).
\end{remark}

\begin{proposition}
(\textbf{Comparison Principle) \label{PrincComp}}Let $W$ be a bounded u.s.c.
viscosity subsolution of (\ref{HJB}) and let $V$ be a bounded l.s.c. viscosity
supersolution of (\ref{HJB}). Moreover, we assume that either $W$ or $V$ is
uniformly continuous. Then%
\[
W\left(  x\right)  \leq V\left(  x\right)  ,
\]
for all $x\in%
\mathbb{R}
^{N}.$
\end{proposition}

The arguments for the Proof are standard. For reader's convenience, we give
the Proof in the Appendix.

The main result of the subsection is the following characterization of the
$\varepsilon$-viability property with respect to the controlled piecewise
deterministic Markov process.

\begin{theorem}
\label{thViab}Given a nonempty, closed set $K\subset%
\mathbb{R}
^{N},$ the following properties are equivalent:

(i) $K$ is $\varepsilon$-viable;

(ii) The following assertions hold simultaneously:

(a) for every $x\in\partial K,$ and every $p\in N_{K}\left(  x\right)  ,$
\[
\inf_{u\in U}\left\{  \left\langle f\left(  x,u\right)  ,p\right\rangle
+\lambda\left(  x,u\right)  Q\left(  x,u,K^{c}\right)  \right\}  \leq0.
\]

(b) for every $x\in\overset{\circ}{K},$%
\[
\inf_{u\in U}\left\{  \lambda\left(  x,u\right)  Q\left(  x,u,K^{c}\right)
\right\}  \leq0.
\]

\end{theorem}

\begin{proof}
We begin by proving that $(ii)\Rightarrow(i).$ We claim that the function
\[
V\left(  x\right)  =\left\{
\begin{array}
[c]{l}%
0\text{, if }x\in K,\\
1\text{, otherwise.}%
\end{array}
\right.
\]
is a viscosity supersolution for (\ref{HJB}). By definition, $V$ is lower
semi-continuous. Obviously, the supersolution condition holds true for all
$x\in%
\mathbb{R}
^{N}\smallsetminus\partial K$. Let us now fix a point $x\in\partial K$. If
$\varphi\in C_{b}^{1}\left(  \mathcal{N}_{x}\right)  ,$ for some
$\mathcal{N}_{x}\subset%
\mathbb{R}
^{N}$ neighborhood of $x,$ is such that $\left(  V-\varphi\right)  $ admits a
global minimum at $x,$ then $\nabla\varphi\left(  x\right)  \in\mathcal{N}%
_{K}\left(  x\right)  $. Thus, the condition (ii) yields
\[
V\left(  x\right)  -\left(  d_{K}\left(  x\right)  \wedge1\right)  +H\left(
x,\nabla\varphi\left(  x\right)  ,V\right)  =-\inf_{u\in U}\left\{
\left\langle f\left(  x,u\right)  ,\nabla\varphi\left(  x\right)
\right\rangle +\lambda\left(  x,u\right)  Q\left(  x,u,K^{c}\right)  \right\}
\geq0.
\]
It follows that $V$ is a bounded viscosity supersolution for (\ref{HJB}).
Using the Comparison Principle, we get
\[
v(x)\leq V(x)=0,
\]
for all $x\in K$ and the $\varepsilon$-viability of $K$ follows.

To prove the converse, we introduce, for every $m\in%
\mathbb{N}
^{\ast},$ the value function $v_{m}$, defined by%
\[
v_{m}(x)=mv\left(  x\right)  =\inf_{u\in\mathbb{L}^{0}\left(
\mathbb{R}
^{N}\times%
\mathbb{R}
_{+};U\right)  }\mathbb{E}\left[  \int_{0}^{\infty}me^{-t}\left(  d_{K}\left(
X_{t}^{x,u}\right)  \wedge1\right)  dt\right]  ,
\]
for all $x\in%
\mathbb{R}
^{N}$. Then, Theorem 1.1 in \cite{S} yields that $v_{m}$ is the unique bounded
viscosity solution of
\begin{equation}
v_{m}\left(  x\right)  -m\left(  d_{K}\left(  x\right)  \wedge1\right)
+H\left(  x,\nabla v_{m}\left(  x\right)  ,v_{m}\right)  =0, \label{HJm}%
\end{equation}
where the Hamiltonian $H$ is given by (\ref{H}).

\underline{Step 1.} We claim that there exists a positive constant $c>0$ such
that, for all $x\in K^{c}$ and all $m\geq1,$
\begin{equation}
v_{m}\left(  x\right)  \geq mc\left(  d_{K}\left(  x\right)  \wedge1\right)
^{2}. \label{ineqvm}%
\end{equation}
We recall that on the set $\left\{  T_{1}>t\right\}  ,$ $X_{t}^{x,u}=\Phi
_{t}^{0,x,u}.$ Standard estimates yield the existence of a positive constant
$c_{1}$ which is independent of $x$, $u$ and $t$ such that
\[
\left\vert \Phi_{t}^{0,x,u}-x\right\vert \leq c_{1}t,
\]
for all $t\geq0.$ Thus, on the set $\left\{  T_{1}>\frac{d_{K}\left(
x\right)  \wedge1}{2c_{1}}\right\}  $ one gets
\[
d_{K}\left(  X_{s}^{x,u}\right)  \wedge1\geq\frac{d_{K}\left(  x\right)
\wedge1}{2}>0,
\]
for all $s\leq\frac{d_{K}\left(  x\right)  \wedge1}{2c_{1}}.$ Using the
Assumptions A1-A2, one easily proves that%
\begin{align*}
\mathbb{E}\left[  \int_{0}^{\infty}me^{-t}\left(  d_{K}\left(  X_{t}%
^{x,u}\right)  \wedge1\right)  dt\right]   &  \geq mE\left[  \int_{0}%
^{\frac{d_{K}\left(  x\right)  \wedge1}{2c_{1}}}e^{-t}\frac{d_{K}\left(
x\right)  \wedge1}{2}dt1_{\left\{  T_{1}>\frac{d_{K}\left(  x\right)  \wedge
1}{2c_{1}}\right\}  }\right] \\
&  \geq Cm\left(  d_{K}\left(  x\right)  \wedge1\right)  ^{2}.
\end{align*}
Hence, (\ref{ineqvm}) holds true for all $x\in K^{c}.$

\underline{Step 2.} Let us fix $x\in\partial K.$ We consider an arbitrary
$p\in N_{K}\left(  x\right)  $ and introduce the test function
\[
\varphi\left(  y\right)  =\left\langle p,y-x\right\rangle -m^{\frac{1}{4}%
}\left\vert y-x\right\vert ^{2},
\]
for all $y\in\mathcal{%
\mathbb{R}
}^{N}.$ We let $x_{m}\in\overline{B}\left(  x,2\right)  $ such that
\begin{equation}
v_{m}\left(  x_{m}\right)  -\varphi\left(  x_{m}\right)  \leq v_{m}\left(
y\right)  -\varphi\left(  y\right)  ,\label{argmin}%
\end{equation}
for all $y\in\overline{B}\left(  x,2\right)  $. One notices that, for large
enough $m,$ $x_{m}\in\overline{B}\left(  x,1\right)  $. Indeed, this is a
simple consequence of the fact that $v_{m}\left(  x\right)  =\varphi\left(
x\right)  =0$ and, thus,
\begin{equation}
0\leq v_{m}\left(  x_{m}\right)  \leq\varphi\left(  x_{m}\right)
\leq\left\langle p,x_{m}-x\right\rangle -m^{\frac{1}{4}}\left\vert
x_{m}-x\right\vert ^{2}.\label{vmphi}%
\end{equation}
Moreover, for large enough $m,$ the inequality (\ref{argmin}) holds true for
all $y\in%
\mathbb{R}
^{N}.$ The inequalities (\ref{ineqvm}) and (\ref{vmphi}) yield%
\begin{equation}
mc\left(  d_{K}\left(  x_{m}\right)  \wedge1\right)  ^{2}\leq v_{m}\left(
x_{m}\right)  \leq\left\langle p,x_{m}-x\right\rangle -m^{\frac{1}{4}%
}\left\vert x_{m}-x\right\vert ^{2}\nonumber
\end{equation}
This implies
\begin{equation}
\lim_{m\rightarrow\infty}m\left(  d_{K}\left(  x_{m}\right)  \wedge1\right)
^{2}=0\text{, }\lim_{m\rightarrow\infty}x_{m}=x\text{ and }\lim_{m\rightarrow
\infty}m^{\frac{1}{4}}\left\vert x_{m}-x\right\vert ^{2}=0,\label{Case1.2}%
\end{equation}
and%
\begin{equation}
\lim_{m\rightarrow\infty}v_{m}\left(  x_{m}\right)  =0.\label{Case1.3}%
\end{equation}
We claim that
\begin{equation}
\underset{m\rightarrow\infty}{\lim\sup}m^{\frac{1}{4}}\left\vert
x_{m}-x\right\vert =0.\label{Case1.3'}%
\end{equation}
We assume that, on the contrary, there exists some positive real constant
$\delta>0$ such that
\begin{equation}
m^{\frac{1}{4}}\left\vert x_{m}-x\right\vert >\delta,\label{Case1.3.2}%
\end{equation}
for all $m\geq1$. For every $m\geq1$, we choose some $y_{m}\in K$ such that
\begin{equation}
d_{K}\left(  x_{m}\right)  =\left\vert x_{m}-y_{m}\right\vert
.\label{Case1.3.3}%
\end{equation}
The equalities (\ref{Case1.2}) imply that $\underset{m\rightarrow\infty}{\lim
}$ $\ y_{m}=x.$ Together with the choice of $p\in N_{K}\left(  x\right)  ,$
the last limit yields%
\begin{equation}
\left\langle p,y_{m}-x\right\rangle \leq\frac{\delta}{2}\left\vert
y_{m}-x\right\vert ,\label{Case1.3.4}%
\end{equation}
for every $m$ large enough. To simplify the notation, we assume that
(\ref{Case1.3.4}) holds true for all $m\geq1$. Using the inequalities
(\ref{vmphi}), (\ref{Case1.3.4}) and (\ref{Case1.3.2}), we have%
\begin{align*}
0 &  \leq\left\langle p,x_{m}-x\right\rangle -m^{\frac{1}{4}}\left\vert
x_{m}-x\right\vert ^{2}\\
&  \leq\left\langle p,y_{m}-x\right\rangle +\left\langle p,x_{m}%
-y_{m}\right\rangle -m^{\frac{1}{4}}\left\vert x_{m}-x\right\vert \left(
\left\vert y_{m}-x\right\vert -\left\vert x_{m}-y_{m}\right\vert \right)  \\
&  \leq\frac{\delta}{2}\left\vert y_{m}-x\right\vert +\left\vert p\right\vert
d_{K}\left(  x_{m}\right)  -\delta\left\vert y_{m}-x\right\vert +m^{\frac
{1}{4}}d_{K}\left(  x_{m}\right)  \left\vert x_{m}-x\right\vert .
\end{align*}
Therefore,%
\begin{align}
\delta &  <m^{\frac{1}{4}}\left\vert x_{m}-x\right\vert \leq m^{\frac{1}{4}%
}\left(  \left\vert y_{m}-x\right\vert +d_{K}\left(  x_{m}\right)  \right)
\nonumber\\
&  \leq\left(  \frac{2}{\delta}\left\vert p\right\vert +1\right)  m^{\frac
{1}{4}}d_{K}\left(  x_{m}\right)  +\frac{2}{\delta}m^{\frac{1}{2}}d_{K}\left(
x_{m}\right)  \left\vert x_{m}-x\right\vert .\label{Case1.3.5}%
\end{align}
We allow $m\rightarrow\infty$ in the inequality (\ref{Case1.3.5}) and recall
that (\ref{Case1.2}) holds true to come to a contradiction. It follows that
(\ref{Case1.3'}) must hold true.

We recall that the function $v_{m}$ is a bounded, continuous viscosity
supersolution of (\ref{HJm}) to get%
\begin{align}
&  v_{m}\left(  x_{m}\right)  -m\left(  d_{K}\left(  x_{m}\right)
\wedge1\right)  \nonumber\\
&  \geq\inf_{u\in U}\left\{  \left\langle p,f\left(  x_{m},u\right)
\right\rangle -2m^{\frac{1}{4}}\left\langle x_{m}-x,f\left(  x_{m},u\right)
\right\rangle +\lambda\left(  x_{m},u\right)  \int_{%
\mathbb{R}
^{N}}\left(  v_{m}\left(  z\right)  -v_{m}\left(  x_{m}\right)  \right)
Q\left(  x_{m},u,dz\right)  \right\}  .\nonumber
\end{align}
Assumption A1 yieds
\begin{align}
&  \left\langle p,f\left(  x,u\right)  \right\rangle \nonumber\\
&  =\left\langle p,f\left(  x_{m},u\right)  \right\rangle -2m^{\frac{1}{4}%
}\left\langle x_{m}-x,f\left(  x_{m},u\right)  \right\rangle \text{
}+\left\langle p,f\left(  x,u\right)  -f\left(  x_{m},u\right)  \right\rangle
+2m^{\frac{1}{4}}\left\langle x_{m}-x,f\left(  x_{m},u\right)  \right\rangle
\nonumber\\
&  \leq\left\langle p,f\left(  x_{m},u\right)  \right\rangle -2m^{\frac{1}{4}%
}\left\langle x_{m}-x,f\left(  x_{m},u\right)  \right\rangle +C\left\vert
p\right\vert \left\vert x_{m}-x\right\vert +Cm^{\frac{1}{4}}\left\vert
x_{m}-x\right\vert ,\label{Case1.5}%
\end{align}
for all $m\geq1$ and all $u\in U$. Here $C$ is a generic real positive
constant that is independent of $m\geq1$ and $u\in U$ and may change from one
line to another. Let us fix $m_{0}\geq1$. Then, for all $m\geq m_{0}$ and all
$u\in U$, by Assumptions A2-A3, we obtain
\begin{align}
&  \lambda\left(  x_{m},u\right)  \int_{%
\mathbb{R}
^{N}}mc\left(  d_{K}\left(  z\right)  \wedge1\right)  ^{2}Q\left(
x_{m},u,dz\right)  -\lambda\left(  x_{m},u\right)  v_{m}\left(  x_{m}\right)
\nonumber\\
&  \geq cm_{0}\lambda\left(  x,u\right)  \int_{%
\mathbb{R}
^{N}}\left(  d_{K}\left(  z\right)  \wedge1\right)  ^{2}Q\left(
x_{m},u,dz\right)  -C\left\vert x_{m}-x\right\vert m_{0}-Cv_{m}\left(
x_{m}\right)  \nonumber\\
&  \geq cm_{0}\lambda\left(  x,u\right)  \int_{%
\mathbb{R}
^{N}}\left(  d_{K}\left(  z\right)  \wedge1\right)  ^{2}Q\left(
x,u,dz\right)  -m_{0}C\eta_{\left(  d_{K}\wedge1\right)  ^{2}}\left(
\left\vert x_{m}-x\right\vert \right)  -C\left\vert x_{m}-x\right\vert
m_{0}-Cv_{m}\left(  x_{m}\right)  \nonumber\\
&  \geq\lambda\left(  x,u\right)  Q\left(  x,u,K^{c}\right)  -CQ\left(
x,u,K_{m_{0}}\right)  -m_{0}C\eta_{\left(  d_{K}\wedge1\right)  ^{2}}\left(
\left\vert x_{m}-x\right\vert \right)  -C\left\vert x_{m}-x\right\vert
m_{0}-Cv_{m}\left(  x_{m}\right)  ,\label{Case1.6}%
\end{align}
where we use the notation
\[
K_{m_{0}}=\left\{  z\in K^{c}:d_{K}\left(  z\right)  <\frac{1}{\sqrt{m_{0}c}%
}\right\}  .
\]
Finally, using (\ref{Case1.5}) and (\ref{Case1.6}), we get%
\begin{align}
&  \left\{  \left\langle p,f\left(  x,u\right)  \right\rangle +\lambda\left(
x,u\right)  Q\left(  x,u,K^{c}\right)  \right\}  \nonumber\\
&  \leq\left\langle p,f\left(  x_{m},u\right)  \right\rangle -2m^{\frac{1}{4}%
}\left\langle x_{m}-x,\text{ }f\left(  x_{m},u\right)  \right\rangle
+C\left\vert p\right\vert \left\vert x_{m}-x\right\vert +Cm^{\frac{1}{4}%
}\left\vert x_{m}-x\right\vert \nonumber\\
&  +\lambda\left(  x_{m},u\right)  \int_{%
\mathbb{R}
^{N}}\left(  v_{m}\left(  z\right)  -v_{m}\left(  x_{m}\right)  \right)
Q\left(  x_{m},u,dz\right)  +C\sup_{u\in U}Q\left(  x,u,K_{m_{0}}\right)
\nonumber\\
&  +m_{0}C\eta_{\left(  d_{K}\wedge1\right)  ^{2}}\left(  \left\vert
x_{m}-x\right\vert \right)  +C\left\vert x_{m}-x\right\vert m_{0}%
+Cv_{m}\left(  x_{m}\right)  ,\label{Case1.7}%
\end{align}
for all $m\geq m_{0}$ and all $u\in U$. We take in (\ref{Case1.7}) the infimum
over $u\in U,$ then $\lim\sup$ as $m\rightarrow\infty$ and recall that the
inequalities (\ref{Case1.2}), (\ref{Case1.3}), (\ref{Case1.3'}) hold true, to
have%
\begin{equation}
\inf_{u\in U}\left\{  \left\langle p,f\left(  x,u\right)  \right\rangle
+\lambda\left(  x,u\right)  Q\left(  x,u,K^{c}\right)  \right\}  \leq
C\sup_{u\in U}Q\left(  x,u,K_{m_{0}}\right)  \label{Case1.8v2}%
\end{equation}
for all $m_{0}\geq1.$ Notice that $\left(  K_{m_{0}}\right)  $ is a decreasing
sequence of sets such that $\cap_{m_{0}\geq1}K_{m_{0}}=\phi.$ Then, using the
Assumption A4, the inequality (\ref{Case1.8v2}) yields
\begin{equation}
\inf_{u\in U}\left\{  \left\langle p,f\left(  x,u\right)  \right\rangle
+\lambda\left(  x,u\right)  Q\left(  x,u,K^{c}\right)  \right\}
\leq0.\label{Cond1}%
\end{equation}

\underline{Step 3.} For $x\in\overset{\circ}{K},$ we take the test function
$\varphi\left(  y\right)  =-\left\vert y-x\right\vert ^{2}$ for all
$y\in\mathcal{%
\mathbb{R}
}^{N}.$ The same arguments as in Step 2 give
\begin{equation}
\inf_{u\in U}\left\{  \lambda\left(  x,u\right)  Q\left(  x,u,K^{c}\right)
\right\}  \leq0. \label{Cond2}%
\end{equation}

\end{proof}

\subsection{Conditions for invariance\label{Subsection1.2}}

Another problem, closely related to viability is the invariance of a nonempty,
closed set $K\subset%
\mathbb{R}
^{N}.$ Whenever this property is satisfied, the controlled PDMP remains in $K$
independently on the control process and as soon as the initial datum $x\in
K$. Suppose that the initial states of the model to which the PDMP is
associated are known. Then, one may reduce the complexity by restricting the
model to some invariant set containing the initial data. We begin by recalling
the notion of invariance.

\begin{definition}
A nonempty, closed set $K\subset%
\mathbb{R}
^{N}$ is said to be (strongly) invariant with respect to the piecewise
deterministic Markov process $X$ if, for every initial point $x\in K$ and
every admissible control process $u$, $X_{t}^{x,u}\in K,$ $\mathbb{P}$-almost
surely, for all $t\in%
\mathbb{R}
_{+}.$
\end{definition}

The invariance property is related to an optimal control problem for which the
value function $v_{inv}$ is given by
\begin{equation}
v_{inv}\left(  x\right)  =\sup_{u\in\mathbb{L}^{0}\left(
\mathbb{R}
^{N}\times%
\mathbb{R}
_{+};U\right)  }\mathbb{E}\left[  \int_{0}^{\infty}e^{-t}\left(  d_{K}\left(
X_{t}^{x,u}\right)  \wedge1\right)  dt\right]  ,\label{vinv}%
\end{equation}
for all $x\in%
\mathbb{R}
^{N}.$ The main result of the section is

\begin{theorem}
\label{thInv}Let $K\subset%
\mathbb{R}
^{N}$ be a nonempty, closed subset. The following statements are equivalent:

(i) The set $K$ is invariant;

(ii) The value function $v_{inv}\left(  x\right)  =0,$ for all $x\in K.$

(iii) The following conditions hold simultaneously:

(a) for every $x\in\partial K,$ every $p\in N_{K}\left(  x\right)  ,$ and
every $u\in U,$
\[
\left\langle f\left(  x,u\right)  ,p\right\rangle +\lambda\left(  x,u\right)
Q\left(  x,u,K^{c}\right)  \leq0.
\]

(b) for every $x\in\overset{\circ}{K},$and every $u\in U,$%
\[
\lambda\left(  x,u\right)  Q\left(  x,u,K^{c}\right)  \leq0.
\]

\end{theorem}

\begin{proof}
We only need to prove that $(ii)$ and $(iii)$ are equivalent. We begin by
proving that $(iii)$ implies $(ii).$ By Theorem 1.1 in \cite{S}, the function
\[
w=-v_{inv}%
\]
is the unique bounded viscosity solution of the Hamilton-Jacobi
integro-differential equation%
\begin{equation}
w\left(  x\right)  +d_{K}\left(  x\right)  \wedge1+H\left(  x,\nabla
w,w\right)  =0, \label{HJinv}%
\end{equation}
where the Hamiltonian is given by (\ref{H}). As in the proof of Theorem
\ref{thViab}, one notices that the function
\[
V\left(  x\right)  =-1_{K^{c}}(x),\text{ for all }x\in%
\mathbb{R}
^{N}%
\]
is a viscosity subsolution of (\ref{HJinv}). By the comparison principle, we
get that
\[
v_{inv}\left(  x\right)  \leq0,
\]
for all $x\in K.$ The statement follows. The proof of the converse relies on
the same arguments as Steps 1-3 of Theorem \ref{thViab}.
\end{proof}

\section{Reachability of open sets\label{Section2}}

Stability issues are very important for biological networks. For deterministic
models, one can easily decide whether the system is stable, bistable, etc.
However, the behavior is much less obvious for a piecewise deterministic
approach. One should expect that the trajectories of the controlled PDMP
starting from some region around the stable point converge to it.
Alternatively, a point for which arbitrarily small surrounding regions are
invariant (or at least viable) is a good candidate for stability. Thus, the
issue of stability may be addressed via viability techniques. In the case of
multiple stable points, given an arbitrary initial state, it would be
interesting to know to which of these regions the trajectories of the PDMP are
directed. The goal of this section is to address the problem of reachability.

Let us consider an arbitrary nonempty, open set $\mathcal{O}\subset%
\mathbb{R}
^{N}.$ As in the case of viability, the techniques we use rely on the theory
of viscosity solutions for a class of Hamilton-Jacobi integro-differential
equations. We are going to introduce a slight difference in our coefficients
allowing to consider a control couple. To this purpose, we make the following
notations: We let the vector field $\widetilde{f}:%
\mathbb{R}
^{N}\times U\times\overline{B}\left(  0,1\right)  \longrightarrow%
\mathbb{R}
^{N}$ be given by
\begin{equation}
\widetilde{f}\left(  x,u^{1},u^{2}\right)  =f\left(  x+u^{2},u^{1}\right)  ,
\label{ftilda}%
\end{equation}
for all $x\in%
\mathbb{R}
^{N}$, $u^{1}\in U$ and $u^{2}\in\overline{B}\left(  0,1\right)  .$ Similarly,
the function $\widetilde{\lambda}:%
\mathbb{R}
^{N}\times U\times\overline{B}\left(  0,1\right)  \longrightarrow%
\mathbb{R}
_{+}$ is given by
\begin{equation}
\widetilde{\lambda}\left(  x,u^{1},u^{2}\right)  =\lambda\left(  x+u^{2}%
,u^{1}\right)  , \label{lambdatilda}%
\end{equation}
and
\[
\widetilde{Q}\left(  x,u^{1},u^{2},A\right)  =Q\left(  x+u^{2},u^{1}%
,A+u^{2}\right)  ,
\]
where $A+u^{2}=\left\{  a+u^{2}:a\in A\right\}  ,$ for all $x\in%
\mathbb{R}
^{N}$, $u^{1}\in U$, $u^{2}\in\overline{B}\left(  0,1\right)  $ and all Borel
set $A\subset%
\mathbb{R}
^{N}.$

\begin{remark}
1. It is obvious that, for every $h\in C_{b}\left(
\mathbb{R}
^{N}\right)  $ and every $x\in%
\mathbb{R}
^{N}$, $u^{1}\in U$, $u^{2}\in\overline{B}\left(  0,1\right)  ,$%
\[
\int_{%
\mathbb{R}
^{N}}h\left(  z\right)  \widetilde{Q}\left(  x,u^{1},u^{2},dz\right)  =\int_{%
\mathbb{R}
^{N}}h\left(  z-u^{2}\right)  Q\left(  x+u^{2},u^{1},dz\right)  .
\]

2. One can easily check that the assumptions (A1)-(A2) and (A5) hold true for
the characteristic $\left(  \widetilde{f},\widetilde{\lambda},\widetilde
{Q}\right)  $ replacing $\left(  f,\lambda,Q\right)  $ and the set of control
$U$ replaced by $U\times\overline{B}\left(  0,1\right)  .$
\end{remark}

Throughout the section we are going to strengthen (A3) and assume

(B) For each bounded uniformly continuous function $h\in BUC\left(
\mathbb{R}
^{N}\right)  ,$ there exists a continuous function $\eta_{h}:%
\mathbb{R}
\longrightarrow%
\mathbb{R}
$ such that $\eta_{h}\left(  0\right)  =0$ and
\begin{equation}
\sup_{u^{1}\in U,u^{2}\in\overline{B}\left(  0,1\right)  }\left\vert \int_{%
\mathbb{R}
^{N}}h\left(  z-u^{2}\right)  Q\left(  x+u^{2},u^{1},dz\right)  -\int_{%
\mathbb{R}
^{N}}h\left(  z-u^{2}\right)  Q\left(  y+u^{2},u,dz\right)  \right\vert
\leq\eta_{h}\left(  \left\vert x-y\right\vert \right)  . \tag{B}\label{B}%
\end{equation}

\begin{remark}
Similarly to Remark \ref{rem1}, one can alternatively assume

(B') For each bounded uniformly continuous function $h\in BUC\left(
\mathbb{R}
^{N}\right)  ,$ there exists a continuous function $\eta_{h}:%
\mathbb{R}
\longrightarrow%
\mathbb{R}
$ such that $\eta_{h}\left(  0\right)  =0$ and
\begin{align*}
\sup_{u^{1}\in U,u^{2}\in\overline{B}\left(  0,1\right)  }  &  \left\{
\begin{array}
[c]{c}%
\lambda\left(  x+u^{2},u^{1}\right)  \int_{%
\mathbb{R}
^{N}}h\left(  z-u^{2}\right)  Q\left(  x+u^{2},u^{1},dz\right) \\
-\lambda\left(  y+u^{2},u^{1}\right)  \int_{%
\mathbb{R}
^{N}}h\left(  z-u^{2}\right)  Q\left(  y+u^{2},u,dz\right)
\end{array}
\right\} \\
&  \leq\eta_{h}\left(  \left\vert x-y\right\vert \right)  .
\end{align*}

\end{remark}

For every $\varepsilon>0,$ we denote by $\mathcal{E}^{\varepsilon}$ the class
of measurable processes $u^{2}:%
\mathbb{R}
^{N}\times%
\mathbb{R}
_{+}\longrightarrow\overline{B}\left(  0,\varepsilon\right)  .$ For every
admissible control couple $\left(  u^{1},u^{2}\right)  \in\mathbb{L}%
^{0}\left(
\mathbb{R}
^{N}\times%
\mathbb{R}
_{+};U\right)  \mathcal{\times E}^{\varepsilon},$ we let $X_{\cdot}%
^{x,u^{1},u^{2}}$ be the piecewise deterministic process associated to the
characteristic $\left(  \widetilde{f},\widetilde{\lambda},\widetilde
{Q}\right)  .$ Obviously, $X_{\cdot}^{x,u^{1},0}$ is associated to $\left(
f,\lambda,Q\right)  .$

\begin{definition}
Given an initial condition $x\in\mathcal{O}^{c}$ (or even $x\in%
\mathbb{R}
^{N}$), the set $\mathcal{O}$ is reachable starting from $x$ if there exists
some admissible control process $u$ such that the set
\[
\left\{  X_{t}^{x,u,0}\in\mathcal{O},\text{ }t\in\left[  0,\infty\right)
\right\}
\]
has positive probability.
\end{definition}

In connection to this property, we define, for every $\varepsilon\geq0,$ the
value function%
\begin{equation}
v^{\varepsilon}(x)=\inf_{u^{1}\in\mathbb{L}^{0}\left(
\mathbb{R}
^{N}\times%
\mathbb{R}
_{+};U\right)  ,u^{2}\in\mathcal{E}^{\varepsilon}}\mathbb{E}\left[  \int
_{0}^{\infty}-e^{-t}\left(  d_{\mathcal{O}^{c}}\left(  X_{t}^{x,u^{1},u^{2}%
}+u_{t}^{2}\right)  \wedge1\right)  dt\right]  , \label{veps}%
\end{equation}
for all $x\in%
\mathbb{R}
^{N}.$

\begin{remark}
It is obvious that, whenever $v^{0}(x)=0,$ the set $\mathcal{O}$ is not
reachable starting from the point $x.$ On the other hand, whenever
$v^{0}(x)<0,$ there exist a constant $\delta>0,$ an admissible control process
$u_{0}\in\mathbb{L}^{0}\left(
\mathbb{R}
^{N}\times%
\mathbb{R}
_{+};U\right)  $ and $T>0$ such that $\mathbb{E}\left[  \int_{0}^{T}%
e^{-t}\left(  d_{\mathcal{O}^{c}}\left(  X_{t}^{x,u_{0},0}\right)
\wedge1\right)  du\right]  >\delta.$ It follows that the set $\left\{
X_{t}^{x,u_{0},0}\in\mathcal{O},\text{ }t\in\left[  0,T\right]  \right\}  $
must have positive probability. Thus, $\mathcal{O}$ is reachable from $x$ if
and only if $v^{0}(x)<0.$
\end{remark}

Theorem 1.1 in Soner \cite{S} yields that $v^{\varepsilon}$ is the unique
bounded viscosity solution of the following Hamilton-Jacobi
integro-differential equation:%
\begin{align}
0  &  =v^{\varepsilon}\left(  x\right)  +\sup_{\left\vert u^{2}\right\vert
\leq\varepsilon}\left\{  d_{\mathcal{O}^{c}}\left(  x+u_{2}\right)
\wedge1+\sup_{u^{1}\in U}\left\{  -\left\langle f\left(  x+u^{2},u^{1}\right)
,\nabla v^{\varepsilon}\left(  x\right)  \right\rangle \right.  \right.
\nonumber\\
&  \left.  \left.  -\lambda\left(  x+u^{2},u^{1}\right)  \int_{%
\mathbb{R}
^{N}}\left(  v^{\varepsilon}\left(  z\right)  -v^{\varepsilon}\left(
x\right)  \right)  \widetilde{Q}\left(  x,u^{1},u^{2},dz\right)  \right\}
\right\}  , \label{HJeps}%
\end{align}
for all $x\in%
\mathbb{R}
^{N}$. For the particular case $\varepsilon=0$, the value function $v^{0}$ is
the unique bounded uniformly continuous viscosity solution of%
\begin{equation}
v^{0}\left(  x\right)  +d_{\mathcal{O}^{c}}\left(  x\right)  \wedge1+H\left(
x,v^{0}\left(  x\right)  ,v^{0}\right)  =0, \label{HJBreach}%
\end{equation}
for all $x\in%
\mathbb{R}
^{N}$, where the Hamiltonian $H$ is given by (\ref{H}).

\begin{remark}
As a consequence of the definition of $\widetilde{Q},$ for every
$\varepsilon>0$ and every $u^{2}\in\overline{B}\left(  0,\varepsilon\right)
$, the function $w\left(  \cdot\right)  =v^{\varepsilon}\left(  \cdot
-u^{2}\right)  $ is \ a viscosity subsolution of (\ref{HJBreach}).
\end{remark}

We get the following convergence theorem

\begin{theorem}
There exists a decreasing function $\eta:%
\mathbb{R}
_{+}\longrightarrow%
\mathbb{R}
_{+}$ that satisfies $\lim_{\varepsilon\rightarrow0}\eta\left(  \varepsilon
\right)  =0$ and such that
\begin{equation}
\sup_{x\in%
\mathbb{R}
^{N}}\left\vert v^{\varepsilon}(x)-v^{0}\left(  x\right)  \right\vert \leq
\eta\left(  \varepsilon\right)  . \label{uc}%
\end{equation}

\end{theorem}

\begin{proof}
We recall that $v^{0}$ is uniformly continuous and let
\begin{equation}
\omega^{0}\left(  r\right)  =\sup\left\{  \left\vert v^{0}(x)-v^{0}%
(y)\right\vert :x,y\in%
\mathbb{R}
^{N},\left\vert x-y\right\vert \leq r\right\}  , \label{omega0}%
\end{equation}
for all $r>0$ be its continuity modulus. Let us fix $x\in%
\mathbb{R}
^{N}$ and $\varepsilon>0.$ We denote by $\Phi_{\cdot}^{t_{0},x_{0},u^{1}%
,u^{2}}$ the flow associated to the vector field $\widetilde{f}.$ Standard
estimates and the assumption (A1) yield the existence of some positive
constant $C>0$ which is independent of $x$ and $\varepsilon>0$ such that
\begin{equation}
\left\vert \Phi_{t}^{0,x,u^{1},u^{2}}-\Phi_{t}^{0,x,u^{1},0}\right\vert \leq
C\varepsilon, \label{est}%
\end{equation}
for all $t\in\left[  0,1\right]  ,$ and all $\left(  u^{1},u^{2}\right)
\in\mathbb{L}^{0}\left(
\mathbb{R}
^{N}\times%
\mathbb{R}
_{+};U\right)  \times\mathcal{E}^{\varepsilon}.$ The constant $C$ is generic
and may change from one line to another. We emphasize that throughout the
proof, $C$ may be chosen independent of $x\in%
\mathbb{R}
^{N},$ $\varepsilon>0$ and of $\left(  u^{1},u^{2}\right)  \in\mathbb{L}%
^{0}\left(
\mathbb{R}
^{N}\times%
\mathbb{R}
_{+};U\right)  \times\mathcal{E}^{\varepsilon}.$ Using the dynamic programming
principle (Soner \cite{S}, Equation (0.8)), for every admissible control
process $u^{1}\in\mathbb{L}^{0}\left(
\mathbb{R}
^{N}\times%
\mathbb{R}
_{+};U\right)  $, the following inequality holds true
\begin{equation}
v^{0}\left(  x\right)  \leq\mathbb{E}\left[  \int_{0}^{T_{1}\wedge1}%
-e^{-t}\left(  d_{\mathcal{O}^{c}}\left(  X_{t}^{x,u^{1},0}\right)
\wedge1\right)  dt+e^{-T_{1}\wedge1}v^{0}\left(  X_{T_{1}\wedge1}^{x,u^{1}%
,0}\right)  \right]  . \label{I0}%
\end{equation}
We consider an arbitrary admissible control couple $\left(  u^{1}%
,u^{2}\right)  \in\mathbb{L}^{0}\left(
\mathbb{R}
^{N}\times%
\mathbb{R}
_{+};U\right)  \times\mathcal{E}^{\varepsilon}$. For simplicity, we introduce
the following notations:
\begin{align*}
\lambda^{1}\left(  t\right)   &  =\lambda\left(  \Phi_{t}^{0,x,u^{1},0}%
,u_{t}^{1}\right)  ,\text{ }\Lambda^{1}\left(  t\right)  =\exp\left(
-\int_{0}^{t}\lambda^{1}\left(  s\right)  ds\right) \\
\lambda^{1,2}\left(  t\right)   &  =\lambda\left(  \Phi_{t}^{0,x,u^{1},u^{2}%
}+u_{t}^{2},u_{t}^{1}\right)  ,\text{ }\Lambda^{1,2}\left(  t\right)
=\exp\left(  -\int_{0}^{t}\lambda^{1,2}\left(  s\right)  ds\right)  ,
\end{align*}
for all $t\geq0$. We denote the right-hand member of the inequality (\ref{I0})
by $I.$ Then, $I$ is explicitly given by
\begin{align*}
I  &  =\int_{0}^{1}\lambda^{1}(t)\Lambda^{1}\left(  t\right)  \int_{0}%
^{t}-e^{-s}\left(  d_{\mathcal{O}^{c}}\left(  \Phi_{s}^{0,x,u^{1},0}\right)
\wedge1\right)  dsdt\\
&  +\int_{0}^{1}\lambda^{1}(t)\Lambda^{1}\left(  t\right)  e^{-t}\int_{%
\mathbb{R}
^{N}}v^{0}\left(  z\right)  Q\left(  \Phi_{t}^{0,x,u^{1},0},u_{t}%
^{1},dz\right)  dt\\
&  +\Lambda^{1}\left(  1\right)  \int_{0}^{1}-e^{-t}\left(  d_{\mathcal{O}%
^{c}}\left(  \Phi_{t}^{0,x,u^{1},0}\right)  \wedge1\right)  dt+\Lambda
^{1}\left(  1\right)  e^{-1}v^{0}\left(  \Phi_{1}^{0,x,u^{1},0}\right) \\
&  =I_{1}+I_{2}+I_{3}+I_{4}.
\end{align*}
Using the inequality (\ref{est}) and the assumption (A2), one gets%
\begin{align}
I_{1}  &  \leq\int_{0}^{1}\lambda^{1,2}(t)\Lambda^{1,2}\left(  t\right)
\int_{0}^{t}-e^{-s}\left(  d_{\mathcal{O}^{c}}\left(  \Phi_{s}^{x,u^{1},u^{2}%
}+u_{s}^{2}\right)  \wedge1\right)  dsdt+C\varepsilon,\label{I1}\\
I_{3}  &  \leq\Lambda^{1,2}\left(  1\right)  \int_{0}^{1}-e^{-t}\left(
d_{\mathcal{O}^{c}}\left(  \Phi_{t}^{x,u^{1},u^{2}}\right)  \wedge1\right)
dt+C\varepsilon. \label{I2}%
\end{align}
For the term $I_{2},$ one has%
\begin{align}
I_{2}  &  \leq\int_{0}^{1}\lambda^{1,2}(t)\Lambda^{1,2}\left(  t\right)
e^{-t}\int_{%
\mathbb{R}
^{N}}v^{0}\left(  z-u_{t}^{2}\right)  Q\left(  \Phi_{t}^{x,u^{1},u^{2}}%
+u_{t}^{2},u_{t}^{1},dz\right)  dt+C\left(  \varepsilon+W_{v^{0}}\left(
C\varepsilon\right)  +\omega^{0}\left(  \varepsilon\right)  \right)
\nonumber\\
&  \leq\int_{0}^{1}\lambda^{1,2}(t)\Lambda^{1,2}\left(  t\right)  e^{-t}\int_{%
\mathbb{R}
^{N}}v^{\varepsilon}\left(  z\right)  \widetilde{Q}\left(  \Phi_{t}%
^{x,u^{1},u^{2}},u_{t}^{1},u_{t}^{2},dz\right)  dt+C\left(  \varepsilon
+W_{v^{0}}\left(  C\varepsilon\right)  +\omega^{0}\left(  \varepsilon\right)
\right) \nonumber\\
&  +\left(  \int_{0}^{1}\lambda^{1,2}(t)\Lambda^{1,2}\left(  t\right)
e^{-t}dt\right)  \sup_{z\in%
\mathbb{R}
^{N}}\left\vert v^{0}(z)-v^{\varepsilon}(z)\right\vert . \label{I3}%
\end{align}
Finally,
\begin{align}
I_{4}  &  \leq\Lambda^{1,2}(1)e^{-1}v^{0}\left(  \Phi_{1}^{x,u^{1},u^{2}%
}\right)  +C\left(  \omega^{0}\left(  C\varepsilon\right)  +\varepsilon\right)
\nonumber\\
&  \leq\Lambda^{1,2}(1)e^{-1}v^{\varepsilon}\left(  \Phi_{1}^{x,u^{1},u^{2}%
}\right)  +C\left(  \omega^{0}\left(  C\varepsilon\right)  +\varepsilon
\right)  +\Lambda^{1,2}(1)e^{-1}\sup_{z}\left(  v^{0}(z)-v^{\varepsilon
}(z)\right)  . \label{I4}%
\end{align}
We substitute (\ref{I1})-(\ref{I4}) in (\ref{I0}). We take the infimum over
the family of $\left(  u^{1},u^{2}\right)  \in\mathbb{L}^{0}\left(
\mathbb{R}
^{N}\times%
\mathbb{R}
_{+};U\right)  \times\mathcal{E}^{\varepsilon}$ and use the dynamic
programming principle to have
\begin{align*}
v^{0}\left(  x\right)   &  \leq v^{\varepsilon}(x)+C\left(  \varepsilon
+W_{v^{0}}\left(  C\varepsilon\right)  +\omega^{0}\left(  C\varepsilon\right)
\right) \\
&  +\left(  \int_{0}^{1}\lambda^{1,2}(t)\Lambda^{1,2}\left(  t\right)
e^{-t}dt+\Lambda^{1,2}\left(  1\right)  e^{-1}\right)  \sup_{z}\left(
v^{0}(z)-v^{\varepsilon}(z)\right)  .
\end{align*}
We notice that
\[
\int_{0}^{1}\lambda^{1,2}(t)\Lambda^{1,2}\left(  t\right)  e^{-t}%
dt+\Lambda^{1,2}\left(  1\right)  e^{-1}=1-\int_{0}^{1}e^{-\int_{0}%
^{t}\widetilde{\lambda}\left(  \Phi_{s}^{x,u^{1},u^{2}},u_{s}^{1},u_{s}%
^{2}\right)  ds}e^{-t}dt\leq1-e^{-\left(  \lambda^{\max}+1\right)  }.
\]
Thus,
\[
v^{0}\left(  x\right)  -v^{\varepsilon}(x)\leq C\left(  \varepsilon+W_{v^{0}%
}\left(  C\varepsilon\right)  +\omega^{0}\left(  C\varepsilon\right)  \right)
+\left(  1-e^{-\left(  \lambda^{\max}+1\right)  }\right)  \sup_{z}\left(
v^{0}(z)-v^{\varepsilon}(z)\right)  .
\]
The conclusion follows by taking the supremum over $x\in%
\mathbb{R}
^{N}$ and recalling that $C$ is independent of $x$ and $\varepsilon>0.$
\end{proof}

We introduce the function $\mu^{\ast}:%
\mathbb{R}
^{N}\longrightarrow%
\mathbb{R}
$ defined by%
\begin{align}
\mu^{\ast}\left(  x\right)   &  =\sup\left\{  \mu\in%
\mathbb{R}
:\text{ }\exists\varphi\in C_{b}^{1}\left(
\mathbb{R}
^{N}\right)  \text{ such that }\forall\left(  y,u\right)  \in%
\mathbb{R}
^{N}\times U,\right.  \nonumber\\
&  \left.  \mu\leq\mathcal{U}^{u}\varphi\left(  y\right)  -\left(
d_{\mathcal{O}^{c}}\left(  y\right)  \wedge1\right)  +\left(  \varphi\left(
x\right)  -\varphi\left(  y\right)  \right)  \right\}  ,\label{miu*}%
\end{align}
where
\begin{equation}
\mathcal{U}^{u}\varphi\left(  y\right)  =\left\langle \nabla\varphi\left(
y\right)  ,f\left(  y,u\right)  \right\rangle +\lambda\left(  y,u\right)
\int_{%
\mathbb{R}
^{N}}\left(  \varphi\left(  z\right)  -\varphi\left(  y\right)  \right)
Q\left(  y,u,dz\right)  ,\label{Uu}%
\end{equation}
for all $y\in%
\mathbb{R}
^{N}$. This function is inspired by the results in \cite{BGQ}. It corresponds
to the dual form of some linearized formulation for the discounted control
problem. In fact, one can interpret the initial problem by using occupational
measures. In a second step, the set of occupational measures can be enlarged
to a set of measures satisfying appropriate conditions. These conditions
involve the infinitesimal generator of the underlying process and can be
interpreted as a classical constraint. Minimizing on this set leads to the
same value function. Duality techniques then allow to give a formulation much
like $\mu^{\ast}$ (but for generators associated to Brownian diffusion
processes). 

The main result of the section gives the equality between the reachability
value function $v^{0}$ and $\mu^{\ast}.$

\begin{theorem}
\label{thReach}For every $x\in%
\mathbb{R}
^{N},$ the equality%
\begin{equation}
v^{0}(x)=\mu^{\ast}(x) \label{v0EGmiu}%
\end{equation}
holds true.
\end{theorem}

\begin{proof}
We begin by proving that
\begin{equation}
v^{0}(x)\geq\mu^{\ast}(x),\label{v0GEmiu*}%
\end{equation}
for all $x\in%
\mathbb{R}
^{N}.$ We fix $x\in%
\mathbb{R}
^{N}$ and $\left(  \mu,\varphi\right)  \in%
\mathbb{R}
\times C_{b}^{1}\left(
\mathbb{R}
^{N}\right)  $ such that
\[
\mu\leq\mathcal{U}^{u}\varphi\left(  y\right)  -\left(  d_{\mathcal{O}^{c}%
}\left(  y\right)  \wedge1\right)  +\left(  \varphi\left(  x\right)
-\varphi\left(  y\right)  \right)  ,
\]
for all $y\in%
\mathbb{R}
^{N},u\in U.$ Then, for every $u\in\mathbb{L}^{0}\left(
\mathbb{R}
^{N}\times%
\mathbb{R}
_{+};U\right)  ,$
\[
\mu\leq\mathcal{U}^{u_{t}}\varphi\left(  X_{t}^{x,u,0}\right)  +\varphi\left(
x\right)  -\varphi\left(  X_{t}^{x,u,0}\right)  -\left(  d_{\mathcal{O}^{c}%
}\left(  X_{t}^{x,u,0}\right)  \wedge1\right)  ,
\]
for all $t\geq0.$ Using It\^{o}'s formula (cf. Theorem 31.3 in \cite{D}), the
last inequality yields
\begin{align}
\mu &  \leq\lim_{T\rightarrow\infty}\mathbb{E}\left[  \int_{0}^{T}%
e^{-t}\left(  \mathcal{U}^{u_{t}}\varphi\left(  X_{t}^{x,u,0}\right)
-\varphi\left(  X_{t}^{x,u,0}\right)  \right)  dt\right]  \nonumber\\
&  +\varphi\left(  x\right)  +\mathbb{E}\left[  \int_{0}^{\infty}%
-e^{-t}\left(  d_{\mathcal{O}^{c}}\left(  X_{t}^{x,u,0}\right)  \wedge
1\right)  dt\right]  \nonumber\\
&  =\lim_{T\rightarrow\infty}e^{-T}\mathbb{E}\left[  \varphi\left(
X_{T}^{x,u,0}\right)  \right]  +\mathbb{E}\left[  \int_{0}^{\infty}%
-e^{-t}\left(  d_{\mathcal{O}^{c}}\left(  X_{t}^{x,u,0}\right)  \wedge
1\right)  dt\right]  \nonumber\\
&  =\mathbb{E}\left[  \int_{0}^{\infty}-e^{-t}\left(  d_{\mathcal{O}^{c}%
}\left(  X_{t}^{x,u,0}\right)  \wedge1\right)  dt\right]  .\label{Ineq2.1}%
\end{align}
for all $T\geq0.$ We recall the definition (\ref{miu*}) of $\mu^{\ast}\left(
x\right)  $ and the inequality (\ref{v0GEmiu*}) follows from (\ref{Ineq2.1}).

In order to complete the proof of the Theorem, we still have to prove that
\begin{equation}
\mu^{\ast}(x)\geq v^{0}(x). \label{miuGEvreach}%
\end{equation}
Let us consider $\left(  \rho_{\varepsilon}\right)  $ a sequence of standard
mollifiers $\rho_{\varepsilon}\left(  y\right)  =\frac{1}{\varepsilon^{N}}%
\rho\left(  \frac{y}{\varepsilon}\right)  ,$ $y\in%
\mathbb{R}
^{N},$ $\varepsilon>0,$ where $\rho\in C^{\infty}\left(
\mathbb{R}
^{N}\right)  $ is a positive function such that
\[
Supp(\rho)\subset\overline{B}\left(  0,1\right)  \text{ and }\int_{%
\mathbb{R}
^{N}}\rho(x)dx=1.
\]
We introduce the functions%
\begin{equation}
V^{\varepsilon}=v^{\varepsilon}\ast\rho_{\varepsilon}, \label{Veps}%
\end{equation}
for all $\varepsilon>0$. We claim that these functions are (viscosity)
subsolutions of (\ref{HJBreach}). The Proof follows the same arguments as
Lemma 2.7 in Barles, Jakobsen \cite{BaJ}. For convenience, we give the Proof
in the Appendix. Using the fact that $V^{\varepsilon}$ is a subsolution of
(\ref{HJBreach}), one gets%
\[
V^{\varepsilon}(x)\leq\mu^{\ast}\left(  x\right)  .
\]
It follows, from (\ref{uc}) that%
\[
\left(  v^{0}\ast\rho_{\varepsilon}\right)  \left(  x\right)  \leq\mu^{\ast
}\left(  x\right)  +\eta\left(  \varepsilon\right)  .
\]
We allow $\varepsilon\rightarrow0$ in the last inequality, and recall that
$v^{0}$ is continuous, to finally get (\ref{miuGEvreach}). The Proof of the
Theorem is now complete.
\end{proof}

The previous result gives the following interesting characterization of the
reachability of the set $\mathcal{O}:$

\begin{criterion}
\label{critReach}Let $x\in%
\mathbb{R}
^{N}$ be an arbitrary initial state. Then \ the controlled piecewise
deterministic Markov process starting from $x$ reaches $\mathcal{O}$ if and
only if there exists $n\in%
\mathbb{N}
^{\ast}$ such that for every $\varphi\in C_{b}^{1}\left(
\mathbb{R}
^{N}\right)  $ there exists $u\in U,y\in%
\mathbb{R}
^{N}$ such that
\begin{equation}
\mathcal{U}^{u}\varphi\left(  y\right)  -d_{\mathcal{O}^{c}}\left(  y\right)
\wedge1+\left(  \varphi\left(  x\right)  -\varphi\left(  y\right)  \right)
<-n^{-1}. \label{crit}%
\end{equation}

\end{criterion}

\section{Biological examples}

\subsection{Biochemical reactions and mathematical
assumptions\label{Subsection3.1}}

We begin by recalling some rudiments on piecewise deterministic Markov
processes associated to gene networks. For further contributions on gene
networks modelling the reader is referred to \cite{CDR}. We suppose that the
biological evolution is given by a family of genes $\mathcal{G=}\left\{
g_{i}:i=1,N\right\}  $ interacting through a finite set of reactions
$\mathcal{R}$. Every reaction $r\in$ $\mathcal{R}$ can be represented as
\[
\alpha_{1}^{r}g_{1}+\alpha_{2}^{r}g_{2}+...+\alpha_{N}^{r}g_{N}\overset{k_{r}%
}{\longrightarrow}\beta_{1}^{r}g_{1}+...+\beta_{N}^{r}g_{N}%
\]
and it specifies that $\alpha_{i}^{r}$ molecules of $i$ type (with $1\leq
i\leq N$) called reactants interact in order to form the products ($\beta
_{i}^{r}$ molecules of $i$ type, with $1\leq i\leq N$). The reaction does not
occur instantaneously and one needs to specify the reaction speed $k_{r}>0.$
Also, the presence of all species is not required ($\alpha_{i}^{r},$
$\beta_{i}^{r}\in%
\mathbb{N}
,$ for all $1\leq i\leq N$). The species are partitioned in two classes called
continuous, respectively discrete component. This partition (for further
considerations, see \cite{CDR}) induces a partition of the reactions. In sum,
we distinguish between reactions contributing to the continuous flow
($\mathcal{C=}\left\{  1,2,...,M_{1}\right\}  $) and jump reactions
($\mathcal{J=}\left\{  M_{1}+1,...,card\left(  \mathcal{R}\right)  \right\}
$). To every reaction $r\in\mathcal{R}$, one associates

1) a stoichiometric column vector $\theta^{r}=\beta^{r}-\alpha^{r}\in%
\mathbb{R}
^{N},$

2) a propensity function $\lambda_{r}:%
\mathbb{R}
^{N}\longrightarrow%
\mathbb{R}
_{+}.$

For a $\mathcal{C}$-type reaction, $\lambda_{r}\left(  x\right)  =k_{r}%
{\textstyle\prod\limits_{i=1}^{N}}
x_{i}^{\alpha_{i}^{r}},$ for all $x\in%
\mathbb{R}
^{N}.$

For a $\mathcal{J}$-type reaction, one should require further regularity as
$x_{i}\rightarrow0$. The jump mechanism will specify that the number of
molecules of type $i$ diminishes by $\alpha_{i}^{r}.$ Therefore, in order to
insure positive components, rather then introducing $\lambda_{r}\left(
x\right)  $ as for continuous reactions, one could consider%
\[
\lambda_{r}\left(  x\right)  =k_{r}%
{\textstyle\prod\limits_{\substack{i=1\\\alpha_{i}^{r}>0}}^{N}}
x_{i}^{\alpha_{i}^{r}}\chi\left(  \frac{x_{i}}{\alpha_{i}^{r}}\right)  ,
\]
for some regular function $\chi$ such that $0\leq\chi\leq1$, $\chi\left(
y\right)  =0$, for $0\leq y\leq1$ and $\chi\left(  y\right)  =1,$ for
$y\geq1+err$ (where $err$ is a positive constant).

The next step consists in the construction of two matrix $M_{1}$ whose columns
are the vectors $\alpha^{r}$, where $r\in\mathcal{C}$, respectively $M_{2}$
whose columns are the vectors $\alpha^{r}$, where $r\in\mathcal{J}$. The flow
is the given by
\[
f\left(  x\right)  =M_{1}\times\left(  \lambda_{1}\left(  x\right)
,\lambda_{2}\left(  x\right)  ,...,\lambda_{M_{1}}\left(  x\right)  \right)
,
\]
the jump intensity
\[
\lambda\left(  x\right)  =%
{\textstyle\sum\limits_{r\in J}}
\lambda_{r}\left(  x\right)
\]
and, whenever $\lambda\left(  x\right)  >0,$ the transition measure $Q$ is
given by
\[
Q\left(  x,dz\right)  =%
{\textstyle\sum\limits_{r\in\mathcal{J}}}
\frac{\lambda_{r}\left(  x\right)  }{\lambda\left(  x\right)  }\delta
_{x+\theta^{r}}\left(  dz\right)  .
\]
One can suppose that all $\lambda_{r}$ are bounded by a reasonable constant
$\lambda^{\max}>0$, by replacing $\lambda_{r}\left(  x\right)  $ by
$\lambda_{r}\left(  x\right)  \wedge\lambda^{\max}.$ Then, it is obvious that
A1 and A2 hold true. If $h\in BUC\left(
\mathbb{R}
^{N}\right)  $,
\begin{align*}
&  \left\vert \lambda\left(  x+e\right)  \int_{%
\mathbb{R}
^{N}}h\left(  z-e\right)  Q\left(  x+e,dz\right)  -\lambda\left(  y+e\right)
\int_{%
\mathbb{R}
^{N}}h\left(  z-e\right)  Q\left(  y+e,dz\right)  \right\vert \\
&  \leq%
{\textstyle\sum\limits_{r\in\mathcal{J}}}
\left\vert \lambda_{r}\left(  x+e\right)  h\left(  x+\theta^{r}\right)
-\lambda_{r}\left(  y+e\right)  h\left(  y+\theta^{r}\right)  \right\vert \leq
c\left(  \left\vert x-y\right\vert +\omega_{h}\left(  \left\vert
x-y\right\vert \right)  \right)  ,
\end{align*}
for all $x,y\in%
\mathbb{R}
^{N},$ where $c$ depends on the Lipschitz constant of $\lambda_{r}$,
$\lambda^{\max}$ and $\left\Vert h\right\Vert _{\infty}$ and $\omega_{h}$ is
the continuity modulus of $h.$ This implies that B' (and a fortiori A3') hold
true. The assumption A4 is a simple consequence of the fact that $Q\left(
x,\cdot\right)  $ is a probability measure for every $x\in%
\mathbb{R}
^{N}.$ Also, one easily notices that
\[
Q\left(  x,%
\mathbb{R}
^{N}\smallsetminus\overline{B}\left(  x,\sup_{r\in\mathcal{J}}\left\vert
\theta^{r}\right\vert \right)  \right)  =0,
\]
which implies A5. It follows that all the assumptions we have made throughout
the paper are naturally satisfied for piecewise deterministic systems
associated to regulatory gene networks.

\subsection{On/Off Model\label{Subsection3.2}}

A two-state model is often employed to describe different situations in the
molecular biology. Usually, the two states describe either the presence or the
absence of some rare molecular specie. Whenever the gene $\gamma$ is inactive
(represented by $\gamma=0)$, the molecule $X$ degrades at a rate $r_{0}$,
whileas, whenever $\gamma$ is active ($\gamma=1$), the molecule $X$ increases
at a rate proportional to some given $r_{1}$.

From the mathematical point of view, the system will be given by a process
$(X(t),\gamma(t))$ on the state space $E=%
\mathbb{R}
\times\{0;1\}$. $\ $The component $X(t)$ follows a differential dynamic
depending on the hidden variable%
\[
\frac{dX}{dt}=\left\{
\begin{array}
[c]{c}%
-r_{0}(X),\text{ }if\text{ }\gamma(t)=0,\\
\text{ }r_{1}(X),\text{ }if\text{ }\gamma(t)=1,
\end{array}
\right.
\]
where $r_{0}(x)\geq0$ is a bounded, Lipschitz-continuous consumption term and
$r_{1}(x)\geq0$ is a bounded, Lipschitz continuous production term. To be more
precise, the PDMP associated to the model has the characteristic $\left(
f,\lambda,Q\right)  $ given by $f_{\gamma}(x)=-r_{0}\left(  x\right)
(1-\gamma)+r_{1}\left(  x\right)  \gamma,$ $\lambda_{\gamma}(x)=\lambda
_{\gamma},$ $Q(\gamma,x;A)=Q((\gamma,x);A)=\delta_{\left(  \left(
1-\gamma\right)  ,x\right)  }(A),$ for all $\gamma\in\left\{  0,1\right\}  ,$
$x\in%
\mathbb{R}
,$ and all $A\subset%
\mathbb{R}
$. The vector field for the $\gamma$ component can be considered to be $0.$
One should expect $0$-consumption whenever $X=0$ and $\gamma=0$ i.e.
$r_{0}(0)=0$, and no production whenever $X=\alpha_{\max}$ (some maximum
level) and $\gamma=1$, i.e. $r_{1}(\alpha_{\max})=0$. The assumptions A1-A5
are obviously satisfied.

\begin{proposition}
\label{cook1}The set $K=\left[  0,\alpha_{\max}\right]  \times\left\{
0,1\right\}  $ is invariant with respect to the PDP associated to the On/Off Model.
\end{proposition}

\begin{proof}
If $x\in\left[  0,\alpha_{\max}\right]  ,$ then, by the definition of $Q$,%
\[
Q\left(  \left(  0,x\right)  ,K^{c}\right)  =Q\left(  \left(  1,x\right)
,K^{c}\right)  =0.
\]
One notices that $N_{\left[  0,\alpha_{\max}\right]  }\left(  0\right)  =%
\mathbb{R}
_{-}$ and $N_{\left[  0,\alpha_{\max}\right]  }\left(  \alpha_{\max}\right)  =%
\mathbb{R}
_{+}.$ For every $p\leq0,$
\[
pf_{0}\left(  0\right)  =-pr_{0}\left(  0\right)  =0\text{ and }pf_{1}\left(
0\right)  \leq0.
\]
For every $p\geq0,$%
\[
pf_{0}\left(  \alpha_{\max}\right)  =-pr_{0}\left(  \alpha_{\max}\right)
\leq0\text{ and }pf_{1}\left(  \alpha_{\max}\right)  =0.
\]
Thus, by applying Theorem \ref{thInv}, one gets the invariance of $K$.
\end{proof}

\begin{remark}
The arguments of the previous Proposition yield that $\left[  a,b\right]
\times\left\{  0,1\right\}  $ is invariant if and only if
\[
r_{0}(a)=r_{1}(b)=0.
\]
Therefore, in order for a point $x_{0}$ to be a stable point, one should find
a sequence $\varepsilon\searrow0$ such that $r_{0}(x_{0}-\varepsilon
)=r_{1}(x_{0}+\varepsilon)=0.$ In particular, a necessary condition is that
$r_{0}(x_{0})=r_{1}(x_{0})=0.$
\end{remark}

We now focus on the model introduced in \cite{CGT} for stochastic gene
expression and its implications on haploinsufficiency. This basic model of
gene expression, product accumulation and product degradation can be given by
the following reaction system:%
\[%
\begin{tabular}
[c]{|l|}\hline
\textbf{G}\\\hline
\end{tabular}
\ \underset{k_{d}}{\overset{k_{a}}{\rightleftarrows}}%
\begin{tabular}
[c]{|l|}\hline
\textbf{G*}\\\hline
\end{tabular}
\ \overset{J_{p}}{\rightarrow}%
\begin{tabular}
[c]{||l||}\hline\hline
\textbf{P}\\\hline\hline
\end{tabular}
\ \overset{k_{p}}{\rightarrow}%
\]
This model considers a gene to switch randomly between inactive state (G) and
active state (G*). The activation (respectively deactivation) rate is denoted
by $k_{a}$ (respectively $k_{d}$)$.$\ When active, each gene expresses a
product (P) at a rate $J_{p}.$ The product is degraded at rate $k_{p}.$ One
can represent this model as a particular case of the On/Off system by
considering
\begin{equation}
r_{0}(x)=k_{p}x,\text{ }r_{1}(x)=J_{p}-k_{p}x,\text{ }\lambda_{0}=k_{a},\text{
}\lambda_{1}=k_{d},\text{ and }\alpha_{\max}=\frac{J_{p}}{k_{p}}.
\label{cookOnOff}%
\end{equation}

The following result is a consequence of Criterion \ref{critReach}:

\begin{proposition}
\label{cook2}For every real constants $a,b$ such that $0<a<b<\alpha_{\max},$
we let $\mathcal{O=}\left(  a,b\right)  \times\left\{  0,1\right\}  .$ Then,
for every $x\in\left(  0,\alpha_{\max}\right)  ,$ the set $\mathcal{O}$ is
reachable with respect to the PDMP associated to Cook's model starting from
$\left(  0,x\right)  .$
\end{proposition}

\begin{proof}
Let us fix $x\in\left(  0,\alpha_{\max}\right)  $. One seeks to apply
Criterion \ref{critReach}. We reason by contradiction and assume, that, for
every $n\in%
\mathbb{N}
^{\ast}$, there exists $\varphi_{n},\psi_{n}\in C_{b}^{1}\left(
\mathbb{R}
\right)  $ such that
\begin{equation}
\left\{
\begin{array}
[c]{c}%
-n^{-1}\leq-\varphi_{n}^{\prime}\left(  y\right)  k_{p}y+k_{a}\psi_{n}\left(
y\right)  -\left(  1+k_{a}\right)  \varphi_{n}\left(  y\right)  -d_{\left[
0,a\right]  \cup\left[  b,\alpha_{\max}\right]  }\left(  y\right)
+\varphi_{n}\left(  x\right)  ,\\
-n^{-1}\leq\psi_{n}^{\prime}\left(  y\right)  \left(  J_{p}-k_{p}y\right)
-\left(  1+k_{d}\right)  \psi_{n}\left(  y\right)  +k_{d}\varphi_{n}\left(
y\right)  -d_{\left[  0,a\right]  \cup\left[  b,\alpha_{\max}\right]  }\left(
y\right)  +\varphi_{n}\left(  x\right)  ,
\end{array}
\right.  \label{r0}%
\end{equation}
for all $y\in\left[  0,1\right]  .$ We multiply the first inequality by
$y^{\frac{1+k_{a}-k_{p}}{k_{p}}}$ and integrate on $\left(  0,z\right]  ,$ for
$z>0,$ to get%
\[
k_{p}z^{\frac{1+k_{a}}{k_{p}}}\varphi_{n}\left(  z\right)  \leq\frac{k_{p}%
}{1+k_{a}}z^{\frac{1+k_{a}}{k_{p}}}\left(  n^{-1}+\varphi_{n}\left(  x\right)
\right)  +k_{a}\int_{0}^{z}y^{\frac{1+k_{a}-k_{p}}{k_{p}}}\psi_{n}\left(
y\right)  dy-\int_{0}^{z}y^{\frac{1+k_{a}-k_{p}}{k_{p}}}d_{\left[  0,a\right]
\cup\left[  b,\alpha_{\max}\right]  }\left(  y\right)  dy,
\]
or again
\begin{equation}
\varphi_{n}\left(  z\right)  \leq\frac{1}{1+k_{a}}\left(  n^{-1}+\varphi
_{n}\left(  x\right)  \right)  +\frac{k_{a}}{k_{p}}\frac{\int_{0}^{z}%
y^{\frac{1+k_{a}-k_{p}}{k_{p}}}\psi_{n}\left(  y\right)  dy}{z^{\frac{1+k_{a}%
}{k_{p}}}}-\frac{1}{k_{p}}\frac{\int_{0}^{z}y^{\frac{1+k_{a}-k_{p}}{k_{p}}%
}d_{\left[  0,a\right]  \cup\left[  b,\alpha_{\max}\right]  }\left(  y\right)
dy}{z^{\frac{1+k_{a}}{k_{p}}}}, \label{r1}%
\end{equation}
for all $z\in\left(  0,\alpha_{\max}\right]  .$ We multiply the second
inequality in (\ref{r0}) by $\left(  J_{p}-k_{p}y\right)  ^{\frac
{1+k_{d}-k_{p}}{k_{p}}}$ and integrate on $\left[  z,\alpha_{\max}\right)  ,$
for $z<\alpha_{\max},$ to get
\begin{align}
\psi_{n}\left(  z\right)   &  \leq\frac{1}{1+k_{d}}\left(  n^{-1}+\varphi
_{n}\left(  x\right)  \right)  +k_{d}\frac{\int_{z}^{\alpha_{\max}}\left(
J_{p}-k_{p}y\right)  ^{\frac{1+k_{d}-k_{p}}{k_{p}}}\varphi_{n}\left(
y\right)  dy}{\left(  J_{p}-k_{p}z\right)  ^{\frac{1+k_{d}}{k_{p}}}%
}\nonumber\\
&  -\frac{\int_{z}^{1}\left(  J_{p}-k_{p}y\right)  ^{\frac{1+k_{d}-k_{p}%
}{k_{p}}}d_{\left[  0,a\right]  \cup\left[  b,\alpha_{\max}\right]  }\left(
y\right)  dy}{\left(  J_{p}-k_{p}z\right)  ^{\frac{1+k_{d}}{k_{p}}}},
\label{r2}%
\end{align}
for all $z\in\left[  0,\alpha_{\max}\right)  .$ We denote by $a_{n}$ the
maximum value of $\varphi_{n}$ on $\left[  0,\alpha_{\max}\right]  .$ It
follows that
\[
\psi_{n}\left(  z\right)  \leq\frac{1}{1+k_{d}}\left(  n^{-1}+\varphi
_{n}\left(  x\right)  \right)  +\frac{k_{d}}{1+k_{d}}a_{n}-\frac{\int
_{z}^{\alpha_{\max}}\left(  J_{p}-k_{p}y\right)  ^{\frac{1+k_{d}-k_{p}}{k_{p}%
}}d_{\left[  0,a\right]  \cup\left[  b,\alpha_{\max}\right]  }\left(
y\right)  dy}{\left(  J_{p}-k_{p}z\right)  ^{\frac{1+k_{d}}{k_{p}}}},
\]
for all $z\in\left[  0,\alpha_{\max}\right)  .$ We substitute the last
inequality in (\ref{r1}), to have%
\begin{equation}
\varphi_{n}\left(  z\right)  \leq\frac{k_{a}+k_{d}+1}{\left(  1+k_{a}\right)
\left(  1+k_{d}\right)  }\left(  n^{-1}+\varphi_{n}\left(  x\right)  \right)
+\frac{k_{a}k_{d}}{\left(  1+k_{a}\right)  \left(  1+k_{d}\right)  }%
a_{n}-f(z), \label{r3}%
\end{equation}
for all $z\in\left(  0,\alpha_{\max}\right)  .$ The function $f$ is defined
by
\[
f\left(  z\right)  =\frac{1}{k_{p}}\frac{\int_{0}^{z}\left(  y^{\frac
{1+k_{a}-k_{p}}{k_{p}}}d_{\left[  0,a\right]  \cup\left[  b,\alpha_{\max
}\right]  }\left(  y\right)  +k_{a}y^{\frac{1+k_{a}-k_{p}}{k_{p}}}\frac
{\int_{y}^{\alpha_{\max}}\left(  J_{p}-k_{p}u\right)  ^{\frac{1+k_{d}-k_{p}%
}{k_{p}}}d_{\left[  0,a\right]  \cup\left[  b,\alpha_{\max}\right]  }\left(
u\right)  du}{\left(  J_{p}-k_{p}y\right)  ^{\frac{1+k_{d}}{k_{p}}}}\right)
dy}{z^{\frac{1+k_{a}}{k_{p}}}},
\]
for all $z\in\left(  0,\alpha_{\max}\right)  .$ One notices that
\[
\lim_{z\rightarrow0+}f(z)=\frac{k_{a}}{\left(  1+k_{a}\right)  J_{p}%
^{\frac{1+k_{d}}{k_{p}}}}\int_{a}^{b}\left(  J_{p}-k_{p}u\right)
^{\frac{1+k_{d}-k_{p}}{k_{p}}}d_{\left[  0,a\right]  \cup\left[
b,\alpha_{\max}\right]  }\left(  u\right)  du>0.
\]
Therefore, $f\left(  z\right)  >\delta,$ for some positive constant $\delta.$
In particular, the inequality (\ref{r3}) written for $z=x$ gives
\[
\varphi_{n}\left(  x\right)  \leq\frac{k_{a}+k_{d}+1}{k_{a}k_{d}}n^{-1}%
+a_{n}-\frac{\left(  1+k_{a}\right)  \left(  1+k_{d}\right)  }{k_{a}k_{d}%
}\delta.
\]
We return to (\ref{r3}) to obtain%
\begin{equation}
\varphi_{n}\left(  z\right)  \leq\frac{k_{a}+k_{d}+1}{k_{a}k_{d}}n^{-1}%
+a_{n}-\frac{\left(  1+k_{a}\right)  \left(  1+k_{d}\right)  }{k_{a}k_{d}%
}\delta,
\end{equation}
for all $z\in\left(  0,\alpha_{\max}\right)  .$ Finally, taking the supremum
over $z\in\left(  0,\alpha_{\max}\right)  ,$ yields
\[
0<\delta\leq\frac{k_{a}+k_{d}+1}{\left(  1+k_{a}\right)  \left(
1+k_{d}\right)  }n^{-1}.
\]
The last inequality fails to hold for large enough $n$. The assertion of our
Proposition follows.
\end{proof}

\begin{remark}
The reachability result is also true when starting from a generic point
$\left(  1,x\right)  $ replacing $\left(  0,x\right)  .$
\end{remark}

We now illustrate (Figure 1) the viability result in Proposition \ref{cook1}
and the reachability properties given by Proposition \ref{cook2}. We use the
classical description of the PDMP associated to Cook's model. The invariant
set is represented in green ($\left[  0,\alpha_{\max}\right]  $) and we
simulate a trajectory starting from a randomly chosen initial value for the
protein. The time horizon is chosen very small ($100$) and the trajectory is
represented in red. The reachable set is given by randomly generated
$a,b\in\left(  0,\alpha_{\max}\right)  $ and is represented by the blue border
lines. Whenever the sample remains in the target set for two consecutive time
steps, the trajectory is represented in blue.%

\[%
{\parbox[b]{4.3855in}{\begin{center}
\fbox{\includegraphics[
height=3.3018in,
width=4.3855in
]%
{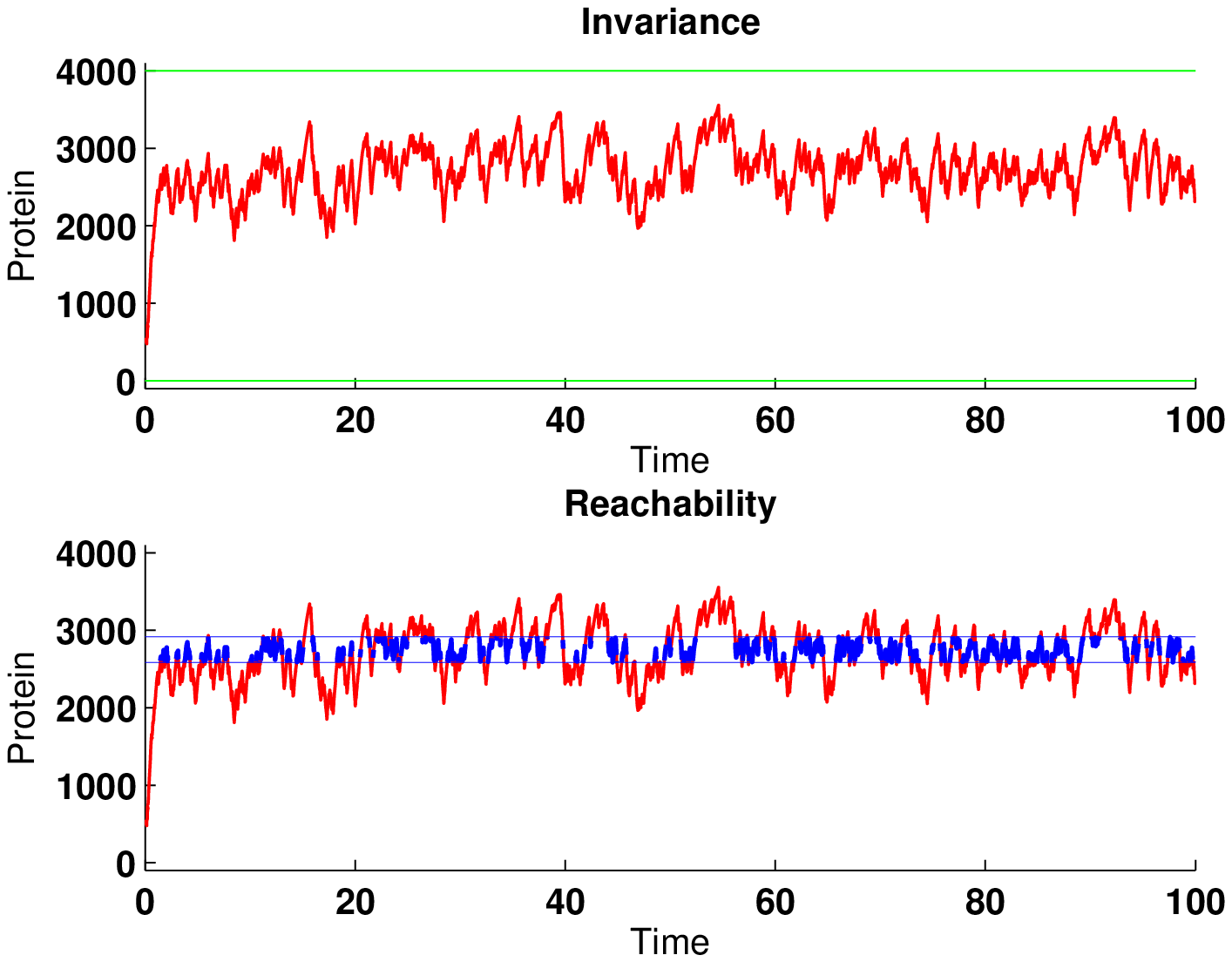}%
}\\
Figure 1. Invariance and reachability properties for Cook's model
\end{center}}}%
\]

\subsection{Bacteriophage $\lambda$\label{Subsection3.3.}}

We consider the model introduced in \cite{HPDC} to describe the regulation of
gene expression. The model is derived from the promoter region of
bacteriophage $\lambda$. The simplification proposed by the authors of
\cite{HPDC} consists in considering a mutant system in which only two operator
sites (known as OR2 and OR3) are present. The gene cI expresses repressor
(CI), which dimerizes and binds to the DNA as a transcription factor in one of
the two available sites. The site OR2 leads to enhanced transcription, while
OR3 represses transcription. Using the notations in \cite{HPDC}, we let $X$
stand for the repressor, $X_{2}$ for the dimer, $D$ for the DNA promoter site,
$DX_{2}$ for the binding to the OR2 site, $DX_{2}^{\ast}$ for the binding to
the OR3 site and $DX_{2}X_{2}$ for the binding to both sites. We also denote
by $P$ the RNA\ polymerase concentration and by $n$ the number of proteins per
mRNA transcript. The dimerization, binding, transcription and degradation
reactions are summarized by%

\[
\left\{
\begin{array}
[c]{l}%
2X\text{ \ \ \ \ \ \ \ \ \ }\overset{K_{1}}{\rightleftarrows}X_{2},\\
D+X_{2}\text{ \ \ \ }\overset{K_{2}}{\rightleftarrows}DX_{2},\\
D+X_{2}\text{ \ \ \ }\overset{K_{3}}{\rightleftarrows}DX_{2}^{\ast},\\
DX_{2}+X_{2}\overset{K_{4}}{\rightleftarrows}DX_{2}X.\\
DX_{2}+P\text{ \ }\overset{K_{t}}{\rightarrow}DX_{2}+P+nX\\
X\text{ \ \ \ \ \ \ \ \ \ \ \ }\overset{K_{d}}{\rightarrow}.
\end{array}
\right.
\]

To this biological system we associate a piecewise deterministic process on
the state space $E=\left\{  \nu\in\left\{  0,1\right\}  ^{4}:\sum_{i=1}^{4}%
\nu_{i}=1\right\}  \times%
\mathbb{R}
^{2}.$ The characteristic is given by
\begin{align*}
f_{\nu}\left(  x_{1},x_{2}\right)   &  =f\left(  x_{1},x_{2}\right)  =\left(
-2k_{1}x_{1}{}^{2}-k_{d}x_{1}+2k_{-1}x_{2},\text{ }k_{1}x_{1}^{2}-k_{-1}%
x_{2}\right)  ,\\
\lambda\left(  \nu,x\right)   &  =k_{2}x_{2}\chi\left(  x_{2}\right)  \nu
_{1}+k_{3}x_{2}\chi\left(  x_{2}\right)  \nu_{1}+k_{4}x_{2}\chi\left(
x_{2}\right)  \nu_{2}+k_{t}\nu_{2}+k_{-2}\nu_{2}+k_{-3}\nu_{3}+k_{-4}\nu
_{4},\\
\lambda\left(  \nu,x\right)  Q\left(  \left(  \lambda,x\right)  ;dz\right)
&  =k_{2}x_{2}\chi\left(  x_{2}\right)  \nu_{1}\delta_{\left(  x_{1}%
,x_{2}-1,\nu_{1}-1,\nu_{2}+1,\nu_{3},\nu_{4}\right)  }\left(  dz\right) \\
&  +k_{3}x_{2}\chi\left(  x_{2}\right)  \nu_{1}\delta_{\left(  x_{1}%
,x_{2}-1,\nu_{1}-1,\nu_{2},\nu_{3}+1,\nu_{4}\right)  }\left(  dz\right) \\
&  +k_{4}x_{2}\chi\left(  x_{2}\right)  \nu_{2}\delta_{\left(  x_{1}%
,x_{2}-1,\nu_{1},\nu_{2}-1,\nu_{3},\nu_{4}+1\right)  }\left(  dz\right) \\
&  +k_{t}\nu_{2}\delta_{\left(  x_{1}+n,x_{2},\nu_{1},\nu_{2},\nu_{3},\nu
_{4}\right)  }\left(  dz\right)  +k_{-2}\nu_{2}\delta_{\left(  x_{1}%
,x_{2}+1,\nu_{1}+1,\nu_{2}-1,\nu_{3},\nu_{4}\right)  }\left(  dz\right) \\
&  +k_{-3}\nu_{3}\delta_{\left(  x_{1},x_{2}+1,\nu_{1}+1,\nu_{2},\nu_{3}%
-1,\nu_{4}\right)  }\left(  dz\right)  +k_{-4}\nu_{4}\delta_{\left(
x_{1},x_{2}+1,\nu_{1},\nu_{2}+1,\nu_{3},\nu_{4}-1\right)  }\left(  dz\right)
,
\end{align*}
for every $\left(  \nu,x\right)  \in E.$ The function $\chi$ is a smooth
function such that $0\leq\chi\leq1$, $\chi\left(  y\right)  =0$ for $y<1$ and
$\chi\left(  y\right)  =1$ for $y\geq1+err.$ We consider the stability
question for this system. Obviously, whenever a point $\left(  \nu^{0}%
,x^{0}\right)  $ is candidate to stability, one should at least expect that
this point should be stable with respect to the deterministic evolution. One
easily notices that the unique equilibrium point for the ordinary equation
driven by the vector field $f$ must satisfy $x^{0}=\left(  0,0\right)  .$
Therefore, any candidate for stability with respect to the piecewise
deterministic evolution associated to the lambda phage should be of this form.

We shall prove that any small enough region surrounding $\left(
0,0,1,0,0,0\right)  $ is invariant with respect to the PDMP. We emphasize that
similar arguments can be used to infer that no other point has similar
stability properties. However, different invariant set may exist. Thus,
bistability of bacteriophage $\lambda$ should be understood as: a stable state
$\left(  0,0,1,0,0,0\right)  $ and some stability (invariance) region.

The main result of the subsection is

\begin{proposition}
For every $\frac{k_{d}^{2}}{4k_{1}k_{-1}}\wedge1>\varepsilon>0,$ the set
\[
K_{\varepsilon}=\left\{  \left(  1,0,0,0\right)  \right\}  \times\left[
0,\frac{2k_{-1}}{k_{d}}\varepsilon\right]  \times\left[  0,\varepsilon\right]
\]
is invariant with respect to the PDMP associated to the bacteriophage
$\lambda$ model.
\end{proposition}

\begin{proof}
We notice that, for every $\left(  x_{1},x_{2}\right)  \in\left[
0,\frac{2k_{-1}}{k_{d}}\varepsilon\right]  \times\left[  0,\varepsilon\right]
,$ one has
\[
\lambda\left(  1,0,0,0,x_{1},x_{2}\right)  =0.
\]
The following table\ gives, for all the possible values of $\left(
x_{1},x_{2}\right)  ,$ the explicit form of the normal cone to $K_{\varepsilon
}$ at $\left(  x_{1},x_{2}\right)  $ and the expression of $\left\langle
p,f\left(  x_{1},x_{2}\right)  \right\rangle $ for every $p=\left(
p_{1},p_{2}\right)  \in\mathbf{N}_{K_{\varepsilon}}\left(  x_{1},x_{2}\right)
.$%

\begin{tabular}
[c]{|l|l|l|l|}\hline
$x_{1}$ & $x_{2}$ & $\mathbf{N}_{K_{\varepsilon}}\left(  x_{1},x_{2}\right)  $
& $\left\langle p,f\left(  x_{1},x_{2}\right)  \right\rangle $\\\hline
$0$ & $0$ & $\left(
\mathbb{R}
_{-}\right)  ^{2}$ & $0$\\\hline
$0$ & $\left(  0,\varepsilon\right)  $ & $%
\mathbb{R}
_{-}\times\left\{  0\right\}  $ & $2p_{1}k_{-1}x_{2}$\\\hline
$0$ & $\varepsilon$ & $%
\mathbb{R}
_{-}\times%
\mathbb{R}
_{+}$ & $\left(  2p_{1}-p_{2}\right)  k_{-1}\varepsilon$\\\hline
$\left(  0,\frac{2k_{-1}}{k_{d}}\varepsilon\right)  $ & $0$ & $\left\{
0\right\}  \times%
\mathbb{R}
_{-}$ & $p_{2}k_{1}x_{1}^{2}$\\\hline
$\left(  0,\frac{2k_{-1}}{k_{d}}\varepsilon\right)  $ & $\left(
0,\varepsilon\right)  $ & $\left\{  0\right\}  \times\left\{  0\right\}  $ &
$0$\\\hline
$\left(  0,\frac{2k_{-1}}{k_{d}}\varepsilon\right)  $ & $\varepsilon$ &
$\left\{  0\right\}  \times%
\mathbb{R}
_{+}$ & $p_{2}\left(  k_{1}x_{1}^{2}-k_{-1}\varepsilon\right)  $\\\hline
$\frac{2k_{-1}}{k_{d}}\varepsilon$ & $0$ & $%
\mathbb{R}
_{+}\times%
\mathbb{R}
_{-}$ & $-p_{1}\left(  \frac{8k_{1}k_{-1}^{2}}{k_{d}^{2}}\varepsilon
^{2}+2k_{-1}\varepsilon\right)  +p_{2}\frac{4k_{1}k_{-1}^{2}}{k_{d}^{2}%
}\varepsilon^{2}$\\\hline
$\frac{2k_{-1}}{k_{d}}\varepsilon$ & $\left(  0,\varepsilon\right)  $ & $%
\mathbb{R}
_{+}\times\left\{  0\right\}  $ & $-p_{1}\left(  \frac{8k_{1}k_{-1}^{2}}%
{k_{d}^{2}}\varepsilon^{2}+2k_{-1}\varepsilon-2k_{-1}x_{2}\right)  $\\\hline
$\frac{2k_{-1}}{k_{d}}\varepsilon$ & $\varepsilon$ & $%
\mathbb{R}
_{+}\times%
\mathbb{R}
_{+}$ & $-p_{1}\frac{8k_{1}k_{-1}^{2}}{k_{d}^{2}}\varepsilon^{2}+p_{2}%
k_{-1}\varepsilon\left(  \frac{4k_{1}k_{-1}}{k_{d}^{2}}\varepsilon-1\right)
$\\\hline
\end{tabular}

The last column allows to conclude that $\left\langle p,f\left(  x_{1}%
,x_{2}\right)  \right\rangle \leq0,$ for all $\left(  x_{1},x_{2}\right)
\in\left[  0,\frac{2k_{-1}}{k_{d}}\varepsilon\right]  \times\left[
0,\varepsilon\right]  $ and all $p\in\mathbf{N}_{K_{\varepsilon}}\left(
x_{1},x_{2}\right)  .$ The conclusion follows from Theorem \ref{thInv}.
\end{proof}

We illustrate the invariance result from the previous Proposition. The reader
is invited to notice that, in the setting of the previous Proposition, the
trajectory should be purely deterministic ($\lambda\left(  1,0,0,0,x_{1}%
,x_{2}\right)  =0,$ for every $\left(  x_{1},x_{2}\right)  \in\left[
0,\frac{2k_{-1}}{k_{d}}\varepsilon\right]  \times\left[  0,\varepsilon\right]
$). We randomly simulate the parameter $\varepsilon<\frac{k_{d}^{2}}%
{4k_{1}k_{-1}}\wedge1$ and a starting point $\left(  x_{1},x_{2}\right)
\in\left[  0,\frac{2k_{-1}}{k_{d}}\varepsilon\right]  \times\left[
0,\varepsilon\right]  .$ The simulated trajectory is represented in red and
the bounds $\frac{2k_{-1}}{k_{d}}\varepsilon$ and $\varepsilon$ are given in
green.%
\[%
{\parbox[b]{4.3855in}{\begin{center}
\fbox{\includegraphics[
height=3.3018in,
width=4.3855in
]%
{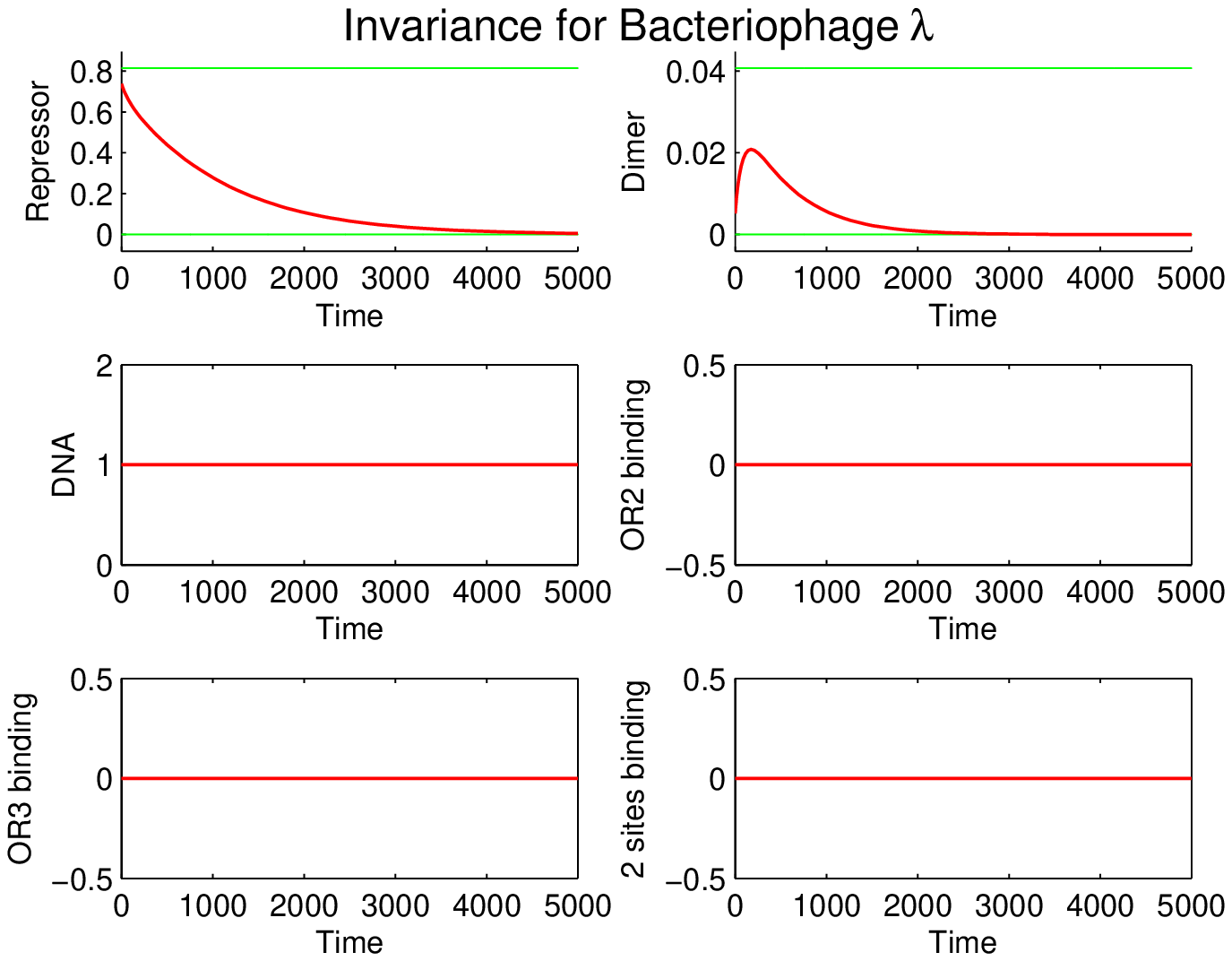}%
}\\
Figure 2. Invariance properties for Bacteriophage $\lambda$
\end{center}}}%
\]

\section{Appendix}

\subsection{A1}

We begin by sketching the proof of the comparison principle.

\begin{proof}
(of Proposition \ref{PrincComp}) Let us suppose that, for some positive
constant $\theta>0,$
\begin{equation}
\sup_{x\in%
\mathbb{R}
^{N}}\left(  W\left(  x\right)  -V\left(  x\right)  \right)  =\theta
>0.\label{uGEv}%
\end{equation}
For every $\delta>0,$ we define%
\[
\theta_{\delta}=\sup_{x\in%
\mathbb{R}
^{N}}\left(  W\left(  x\right)  -V(x)-\delta\left\vert x\right\vert
^{2}\right)  .
\]
It is obvious that $\left(  \theta_{\delta}\right)  $ is decreasing and
$\lim_{\delta\rightarrow0}\theta_{\delta}=\theta.$ For every $\varepsilon>0,$
we introduce
\begin{equation}
\Phi_{\varepsilon,\delta}\left(  x,y\right)  =W\left(  x\right)  -V\left(
y\right)  -\left\vert \frac{x-y}{\varepsilon}\right\vert ^{2}-\delta\left\vert
x\right\vert ^{2},\label{PHIeps}%
\end{equation}
for all $x,y\in%
\mathbb{R}
^{N}.$ We recall that $W$ and $V$ are bounded. Then, using the u.s.c. of
$\Phi_{\varepsilon,\delta},$ we get the existence of some global maximum point
$\left(  x_{\varepsilon,\delta},y_{\varepsilon,\delta}\right)  \in%
\mathbb{R}
^{2N}$ of $\Phi_{\varepsilon,\delta}$. Standard arguments yield the existence
of some $x_{\delta}\in%
\mathbb{R}
^{N}$ such that $\theta_{\delta}=W\left(  x_{\delta}\right)  -V\left(
x_{\delta}\right)  -\delta\left\vert x_{\delta}\right\vert ^{2},$ and
\begin{equation}
\left\{
\begin{array}
[c]{c}%
(i)\text{ }\lim_{\varepsilon\rightarrow0}\left\vert \frac{x_{\varepsilon
,\delta}-y_{\varepsilon,\delta}}{\varepsilon}\right\vert ^{2}=0,\text{
}(ii)\text{ }\lim_{\varepsilon\rightarrow0}x_{\varepsilon,\delta}%
=\lim_{\varepsilon\rightarrow0}y_{\varepsilon,\delta}=x_{\delta},\\
(iii)\text{ }\lim_{\varepsilon\rightarrow0}W\left(  x_{\varepsilon,\delta
}\right)  =W\left(  x_{\delta}\right)  ,\text{ }(iv)\text{ }\lim
_{\varepsilon\rightarrow0}V\left(  y_{\varepsilon,\delta}\right)  =V\left(
x_{\delta}\right)  .
\end{array}
\right.  \label{A1.0}%
\end{equation}
We also obtain%
\begin{equation}
\lim_{\delta\rightarrow0}\delta\left\vert x_{\delta}\right\vert ^{2}%
=0.\label{A1.1}%
\end{equation}
We recall that $W$ is a viscosity subsolution and consider the test function
\[
\varphi_{\varepsilon,\delta}(x)=V\left(  y_{\varepsilon,\delta}\right)
+\left\vert \frac{x-y_{\varepsilon,\delta}}{\varepsilon}\right\vert
^{2}+\delta\left\vert x\right\vert ^{2},
\]
for all $x\in%
\mathbb{R}
^{N}.$ We get%
\begin{align*}
0 &  \geq W\left(  x_{\varepsilon,\delta}\right)  -d_{K}\left(  x_{\varepsilon
,\delta}\right)  \wedge1+\sup_{u\in U}\left\{  -\frac{2}{\varepsilon^{2}%
}\left\langle f\left(  x_{\varepsilon,\delta},u\right)  ,x_{\varepsilon
,\delta}-y_{\varepsilon,\delta}\right\rangle \right.  \\
&  \left.  \text{ \ \ }-2\delta\left\langle f\left(  x_{\varepsilon,\delta
},u\right)  ,x_{\varepsilon,\delta}\right\rangle -\lambda\left(
x_{\varepsilon,\delta},u\right)  \int_{%
\mathbb{R}
^{N}}\left(  W\left(  z\right)  -W\left(  x_{\varepsilon,\delta}\right)
\right)  Q\left(  x_{\varepsilon,\delta},u,dz\right)  \right\}  .\\
&  \geq W\left(  x_{\varepsilon,\delta}\right)  -d_{K}\left(  x_{\varepsilon
,\delta}\right)  \wedge1-C\delta\left\vert x_{\varepsilon,\delta}\right\vert
+\sup_{u\in U}\left\{  -\frac{2}{\varepsilon^{2}}\left\langle f\left(
x_{\varepsilon,\delta},u\right)  ,x_{\varepsilon,\delta}-y_{\varepsilon
,\delta}\right\rangle \right.  \\
&  \left.  \text{ \ \ }-\lambda\left(  x_{\varepsilon,\delta},u\right)  \int_{%
\mathbb{R}
^{N}}\left(  W\left(  z\right)  -W\left(  x_{\varepsilon,\delta}\right)
\right)  Q\left(  x_{\varepsilon,\delta},u,dz\right)  \right\}  ,
\end{align*}
where $C$ is a generic real constant that may change from one line to another.
Standard estimates yield%
\begin{align}
0 &  \geq W\left(  x_{\varepsilon,\delta}\right)  -d_{K}\left(  x_{\varepsilon
,\delta}\right)  \wedge1-C\left(  \delta\left\vert x_{\varepsilon,\delta
}\right\vert +\left\vert W\left(  x_{\varepsilon,\delta}\right)  -W\left(
x_{\delta}\right)  \right\vert +\left\vert x_{\varepsilon,\delta}-x_{\delta
}\right\vert \right)  \nonumber\\
&  \text{ \ \ \ }+\sup_{u\in U}\left\{  -\frac{2}{\varepsilon^{2}}\left\langle
f\left(  x_{\varepsilon,\delta},u\right)  ,x_{\varepsilon,\delta
}-y_{\varepsilon,\delta}\right\rangle \right.  \nonumber\\
&  \left.  \text{ \ \ }-\lambda\left(  x_{\delta},u\right)  \int_{%
\mathbb{R}
^{N}}\left(  W\left(  z\right)  -W\left(  x_{\delta}\right)  \right)  Q\left(
x_{\varepsilon,\delta},u,dz\right)  \right\}  .\label{sub}%
\end{align}
In a similar way, one has
\begin{align}
0 &  \leq V\left(  y_{\varepsilon,\delta}\right)  -d_{K}\left(  y_{\varepsilon
,\delta}\right)  \wedge1+C\left(  \left\vert V\left(  y_{\varepsilon,\delta
}\right)  -V\left(  x_{\delta}\right)  \right\vert +\left\vert y_{\varepsilon
,\delta}-x_{\delta}\right\vert \right)  \nonumber\\
&  +\sup_{u\in U}\left\{  -\frac{2}{\varepsilon^{2}}\left\langle f\left(
y_{\varepsilon,\delta},u\right)  ,x_{\varepsilon,\delta}-y_{\varepsilon
,\delta}\right\rangle -\lambda\left(  x_{\delta},u\right)  \int_{%
\mathbb{R}
^{N}}\left(  V\left(  z\right)  -V\left(  x_{\delta}\right)  \right)  Q\left(
y_{\varepsilon,\delta},u,dz\right)  \right\}  .\label{super2}%
\end{align}
Combining (\ref{sub}) and (\ref{super2}), we get
\begin{align}
0 &  \leq V\left(  y_{\varepsilon,\delta}\right)  -W\left(  x_{\varepsilon
,\delta}\right)  +C\left(  \left\vert \frac{x_{\varepsilon,\delta
}-y_{\varepsilon,\delta}}{\varepsilon}\right\vert ^{2}+\left\vert
x_{\varepsilon,\delta}-x_{\delta}\right\vert +\left\vert y_{\varepsilon
,\delta}-x_{\delta}\right\vert +\right)  \label{A1.2}\\
&  +C\left(  \delta\left\vert x_{\varepsilon,\delta}\right\vert +\left\vert
V\left(  y_{\varepsilon,\delta}\right)  -V\left(  x_{\delta}\right)
\right\vert +\left\vert W\left(  x_{\varepsilon,\delta}\right)  -W\left(
x_{\delta}\right)  \right\vert \right)  \nonumber\\
&  +\sup_{u\in U}\left\{  \lambda\left(  x_{\delta},u\right)  \left(  \int_{%
\mathbb{R}
^{N}}\left(  W\left(  z\right)  -W\left(  x_{\delta}\right)  \right)  Q\left(
x_{\varepsilon,\delta},u,dz\right)  -\int_{%
\mathbb{R}
^{N}}\left(  V\left(  z\right)  -V\left(  x_{\delta}\right)  \right)  Q\left(
y_{\varepsilon,\delta},u,dz\right)  \right)  \right\}  .\nonumber
\end{align}
On the other hand, we notice that%
\begin{align}
&  \int_{%
\mathbb{R}
^{N}}\left(  W\left(  z\right)  -W\left(  x_{\delta}\right)  \right)  Q\left(
x_{\varepsilon,\delta},u,dz\right)  -\int_{%
\mathbb{R}
^{N}}\left(  V\left(  z\right)  -V\left(  x_{\delta}\right)  \right)  Q\left(
y_{\varepsilon,\delta},u,dz\right)  \nonumber\\
&  \leq\int_{%
\mathbb{R}
^{N}}\left(  W\left(  z\right)  -V\left(  z\right)  -W\left(  x_{\delta
}\right)  +V\left(  x_{\delta}\right)  \right)  Q\left(  x_{\varepsilon
,\delta},u,dz\right)  \nonumber\\
&  +\int_{%
\mathbb{R}
^{N}}V\left(  z\right)  \left(  Q\left(  x_{\varepsilon,\delta},u,dz\right)
-Q\left(  y_{\varepsilon,\delta},u,dz\right)  \right)  .\nonumber
\end{align}
Thus, whenever $V$ is continuous,
\begin{align}
&  \int_{%
\mathbb{R}
^{N}}\left(  W\left(  z\right)  -W\left(  x_{\delta}\right)  \right)  Q\left(
x_{\varepsilon,\delta},u,dz\right)  -\int_{%
\mathbb{R}
^{N}}\left(  V\left(  z\right)  -V\left(  x_{\delta}\right)  \right)  Q\left(
y_{\varepsilon,\delta},u,dz\right)  \nonumber\\
&  \leq\int_{%
\mathbb{R}
^{N}}\left(  \delta\left\vert z\right\vert ^{2}\wedge C\right)  Q\left(
x_{\delta},u,dz\right)  +\eta_{\left(  \delta\left\vert \cdot\right\vert
^{2}\wedge C\right)  }\left(  \left\vert x_{\varepsilon,\delta}-x_{\delta
}\right\vert \right)  +\eta_{V}\left(  \left\vert x_{\varepsilon,\delta
}-y_{\varepsilon,\delta}\right\vert \right)  ,\label{A1.3}%
\end{align}
where $\eta_{\delta}$ and $\eta_{V}$ are given by Assumption A3 and are
independent of $u\in U$. Similar estimates hold true if $W$ is continuous. We
substitute (\ref{A1.3}) in (\ref{A1.2}) and take $\lim\sup$ as $\varepsilon
\rightarrow0$ in (\ref{A1.2}) to obtain%
\begin{align*}
0 &  \leq-\theta_{\delta}+\delta\left\vert x_{\delta}\right\vert +C\sup_{u\in
U}\int_{%
\mathbb{R}
^{N}}\left(  \delta\left\vert z\right\vert ^{2}\wedge1\right)  Q\left(
x_{\delta},u,dz\right)  \\
&  \leq-\theta_{\delta}+\delta\left\vert x_{\delta}\right\vert +C\sup_{u\in
U}\left\{  \int_{\overline{B}\left(  x_{\delta},\delta^{-\frac{1}{4}}\right)
}\delta\left\vert z\right\vert ^{2}Q\left(  x_{\delta},u,dz\right)  +Q\left(
x_{\delta},u,%
\mathbb{R}
^{N}\smallsetminus\overline{B}\left(  x_{\delta},\delta^{-\frac{1}{4}}\right)
\right)  \right\}  \\
&  \leq-\theta_{\delta}+C\left(  \delta\left\vert x_{\delta}\right\vert
^{2}+\delta\left\vert x_{\delta}\right\vert +\delta^{\frac{1}{2}}\right)
+\sup_{\left(  x,u\right)  \in%
\mathbb{R}
^{N}\times U}Q\left(  x,u,%
\mathbb{R}
^{N}\smallsetminus\overline{B}\left(  x,\delta^{-\frac{1}{4}}\right)  \right)
.
\end{align*}
We allow $\delta\rightarrow0$ in the last inequality and recall that A5 holds
true to have
\[
0\leq-\theta.
\]
This comes in contradiction with (\ref{uGEv}). The proof is now complete.
\end{proof}

\subsection{A2}

The Proof of Theorem \ref{thReach} relies on the fact that the functions
$V^{\varepsilon}$ defined by (\ref{Veps}) are viscosity subsolutions of the
Hamilton-Jacobi integro-differential equation (\ref{HJBreach}). The proof
adapts the arguments used in Barles, Jakobsen \cite{BaJ} Lemma 2.7. Following
the proof of this Lemma, we introduce, for every $h>0,$ $u^{2}\in%
\mathbb{R}
^{N},$ $Q_{h}^{u^{2}}=u^{2}+\left[  -\frac{h}{2},\frac{h}{2}\right)  ^{N},$
$\rho_{\varepsilon}^{h,u^{2}}=\int_{Q_{h}^{u^{2}}}\rho_{\varepsilon}(y)dy$,
and $I_{h}\left(  x\right)  =\sum_{u^{2}\in h%
\mathbb{Z}
^{N}}\rho_{\varepsilon}^{h,u^{2}}v^{\varepsilon}\left(  x-u^{2}\right)  .$
Thus, $I_{h}$ is a convex combination of bounded, uniformly continuous
viscosity subsolutions of (\ref{HJBreach}). Moreover, by classical results,
the discretization $I_{h}$ converges uniformly to $V^{\varepsilon}.$ To
conclude, we show that viscosity subsolutions are preserved by convex
combination and uniform convergence.

\begin{proposition}
Given two bounded, uniformly continuous viscosity subsolutions $v_{1}$ and
$v_{2}$ of the Equation (\ref{HJBreach}) and two nonnegative real constants
$\lambda_{1},\lambda_{2}\in%
\mathbb{R}
_{+}$ such that $\lambda_{1}+\lambda_{2}=1,$ the convex combination
$\lambda_{1}v_{1}+\lambda_{2}v_{2}$ is still a viscosity subsolution of
(\ref{HJBreach}).
\end{proposition}

\begin{proof}
The assertion is trivial when either $\lambda_{1}=0$ or $\lambda_{2}=0.$ If
$\lambda_{1}\lambda_{2}\neq0,$ we let $\overline{x}\in%
\mathbb{R}
^{N}$ and $\varphi\in C_{b}^{1}\left(  \mathcal{N}_{\overline{x}}\right)  $ be
a test function such that
\begin{equation}
\lambda_{1}v_{1}\left(  \overline{x}\right)  +\lambda_{2}v_{2}\left(
\overline{x}\right)  -\varphi\left(  \overline{x}\right)  \geq\lambda_{1}%
v_{1}\left(  y\right)  +\lambda_{2}v_{2}\left(  y\right)  -\varphi\left(
y\right)  , \label{A2.0.1}%
\end{equation}
for all $y\in%
\mathbb{R}
^{N}.$ We may assume, without loss of generality that $\varphi\in C_{b}\left(
%
\mathbb{R}
^{N}\right)  .$ Indeed, whenever $\varphi$ does not satisfy this assumption,
one can replace it with some $\varphi^{0}$ defined as follows : First, notice
that there exists some $r>0$ such that $B\left(  \overline{x},2r\right)
\subset\mathcal{N}_{\overline{x}}$. We define%
\begin{align*}
\varphi^{0}\left(  y\right)   &  =\left(  \varphi\left(  y\right)
+\lambda_{1}v_{1}\left(  \overline{x}\right)  +\lambda_{2}v_{2}\left(
\overline{x}\right)  -\varphi\left(  \overline{x}\right)  \right)  \chi\left(
y\right) \\
&  +\left(  \lambda_{1}v_{1}\left(  y\right)  +\lambda_{2}v_{2}\left(
y\right)  \right)  \left(  1-\chi\left(  y\right)  \right)  ,
\end{align*}
for all $y\in%
\mathbb{R}
^{N},$ where $\chi$ is a smooth function such that $0\leq\chi\leq1,$
$\chi(y)=1$, if $y\in B\left(  \overline{x},r\right)  $ and $\chi(y)=0$, if
$y\in%
\mathbb{R}
^{N}\smallsetminus B\left(  \overline{x},2r\right)  $. Then (\ref{A2.0.1})
holds true with $\varphi^{0}$ instead of $\varphi.$ The new function
$\varphi^{0}$ also satisfies
\[
\nabla\varphi^{0}\left(  \overline{x}\right)  =\nabla\varphi\left(
\overline{x}\right)  .
\]
We introduce, for every $\varepsilon>0$
\[
\Phi_{\varepsilon}\left(  x,y\right)  =\lambda_{1}v_{1}\left(  x\right)
+\lambda_{2}v_{2}\left(  y\right)  -\lambda_{1}\varphi\left(  x\right)
-\lambda_{2}\varphi\left(  y\right)  -\frac{1}{\varepsilon^{2}}\left\vert
x-y\right\vert ^{2}-\left\vert x-\overline{x}\right\vert ^{2},
\]
for all $x,y\in%
\mathbb{R}
^{N}.$ We recall that the functions $v_{1},v_{2}$ and $\varphi$ are bounded
and continuous. This yields the existence of a global maximum $\left(
x_{\varepsilon},y_{\varepsilon}\right)  $ of $\Phi_{\varepsilon}$. Moreover,
by standard arguments,
\begin{equation}
\lim_{\varepsilon\rightarrow0}x_{\varepsilon}=\lim_{\varepsilon\rightarrow
0}y_{\varepsilon}=\overline{x},\text{ }\lim_{\varepsilon\rightarrow
0}\left\vert \frac{x_{\varepsilon}-y_{\varepsilon}}{\varepsilon}\right\vert
^{2}=0. \label{A2.0.2}%
\end{equation}
We consider the test function $\psi$ given by
\[
\psi\left(  x\right)  =-\lambda_{2}\lambda_{1}^{-1}v_{2}\left(  y_{\varepsilon
}\right)  +\varphi\left(  x\right)  +\lambda_{2}\lambda_{1}^{-1}\varphi\left(
y_{\varepsilon}\right)  +\frac{\lambda_{1}^{-1}}{\varepsilon^{2}}\left\vert
x-y_{\varepsilon}\right\vert ^{2}+\lambda_{1}^{-1}\left\vert x-\overline
{x}\right\vert ^{2},
\]
for all $x\in%
\mathbb{R}
^{N}.$ We recall that the function $v_{1}$ is a viscosity subsolution for
(\ref{HJBreach}). Then,
\[
v_{1}\left(  x_{\varepsilon}\right)  +d_{\mathcal{O}^{c}}\left(
x_{\varepsilon}\right)  \wedge1+H\left(  x_{\varepsilon},\nabla\varphi\left(
x_{\varepsilon}\right)  +\frac{2\lambda_{1}^{-1}}{\varepsilon^{2}}\left(
x_{\varepsilon}-y_{\varepsilon}\right)  +2\lambda_{1}^{-1}\left(
x_{\varepsilon}-\overline{x}\right)  ,v_{1}\right)  \leq0.
\]
Standard estimates yield%
\begin{align}
0  &  \geq v_{1}\left(  \overline{x}\right)  +d_{\mathcal{O}^{c}}\left(
\overline{x}\right)  \wedge1+\sup_{u\in U}\left\{  -\left\langle f\left(
\overline{x},u\right)  ,\nabla\varphi\left(  \overline{x}\right)
\right\rangle -\frac{2\lambda_{1}^{-1}}{\varepsilon^{2}}\left\langle
x_{\varepsilon}-y_{\varepsilon},f\left(  x_{\varepsilon},u\right)
\right\rangle \right. \nonumber\\
&  \left.  -\lambda\left(  \overline{x},u\right)  \int_{%
\mathbb{R}
^{N}}\left(  v_{1}\left(  z\right)  -v_{1}\left(  \overline{x}\right)
\right)  Q\left(  \overline{x},u,dz\right)  \right\} \nonumber\\
&  -C\left(  \left\vert x_{\varepsilon}-\overline{x}\right\vert +\left\vert
v_{1}\left(  x_{\varepsilon}\right)  -v_{1}\left(  \overline{x}\right)
\right\vert +\left\vert \nabla\varphi\left(  x_{\varepsilon}\right)
-\nabla\varphi\left(  \overline{x}\right)  \right\vert +\eta_{v_{1}}\left(
\left\vert x_{\varepsilon}-\overline{x}\right\vert \right)  \right)  .
\label{v1}%
\end{align}
In a similar way, we get
\begin{align}
0  &  \geq v_{2}\left(  \overline{x}\right)  +d_{\mathcal{O}^{c}}\left(
\overline{x}\right)  \wedge1+\sup_{u\in U}\left\{  -\left\langle f\left(
\overline{x},u\right)  ,\nabla\varphi\left(  \overline{x}\right)
\right\rangle +\frac{2\lambda_{2}^{-1}}{\varepsilon^{2}}\left\langle
x_{\varepsilon}-y_{\varepsilon},f\left(  y_{\varepsilon},u\right)
\right\rangle \right. \nonumber\\
&  \left.  -\lambda\left(  \overline{x},u\right)  \int_{%
\mathbb{R}
^{N}}\left(  v_{2}\left(  z\right)  -v_{2}\left(  \overline{x}\right)
\right)  Q\left(  \overline{x},u,dz\right)  \right\} \nonumber\\
&  -C\left(  \left\vert y_{\varepsilon}-\overline{x}\right\vert +\left\vert
v_{2}\left(  y_{\varepsilon}\right)  -v_{2}\left(  \overline{x}\right)
\right\vert +\left\vert \nabla\varphi\left(  y_{\varepsilon}\right)
-\nabla\varphi\left(  \overline{x}\right)  \right\vert +\eta_{v_{2}}\left(
\left\vert y_{\varepsilon}-\overline{x}\right\vert \right)  \right)
\label{v2}%
\end{align}
Finally, using (\ref{v1}), (\ref{v2}) and (\ref{A2.0.2}), and passing to the
limit as $\varepsilon\rightarrow0$, yields
\[
\left(  \lambda_{1}v_{1}+\lambda_{2}v_{2}\right)  \left(  \overline{x}\right)
+d_{\mathcal{O}^{c}}\left(  \overline{x}\right)  \wedge1+H\left(  \overline
{x},\nabla\varphi\left(  \overline{x}\right)  ,\lambda_{1}v_{1}+\lambda
_{2}v_{2}\right)  \leq0.
\]

\end{proof}

These arguments allow to obtain, by recurrence, that any convex combination of
continuous, bounded viscosity subsolutions is still a subsolution for
(\ref{HJBreach}).

\begin{proposition}
(Stability)

Let $\left(  v_{n}\right)  _{n}$ be a sequence of continuous, uniformly
bounded viscosity subsolutions of (\ref{HJBreach}). Moreover, we suppose that
$v_{n}$ converges uniformly on compact sets to some continuous, bounded
function $v$. Then the function $v$ is a viscosity subsolution of
(\ref{HJBreach}).
\end{proposition}

\begin{proof}
We let $x\in%
\mathbb{R}
^{N}$ and $\varphi\in C_{b}^{1}\left(  \mathcal{N}_{x}\right)  $ be a test
function such that $v-\varphi$ has a global maximum at $x.$ As in the previous
proposition, one can assume, without loss of generality, that $\varphi\in
C_{b}\left(
\mathbb{R}
^{N}\right)  .$ Classical arguments yield the existence of some point
$x_{n}\in%
\mathbb{R}
^{N}$ such that
\[
v_{n}\left(  x_{n}\right)  -\varphi\left(  x_{n}\right)  -\left\vert
x_{n}-x\right\vert ^{2}\geq v_{n}\left(  y\right)  -\varphi\left(  y\right)
-\left\vert y-x\right\vert ^{2},
\]
for all $y\in%
\mathbb{R}
^{N}$ and
\[
\lim_{n\rightarrow\infty}x_{n}=x.
\]
We assume, without loss of generality, that $\left\vert x_{n}-x\right\vert
\leq1,$ and $x_{n}\in\mathcal{N}_{x},$ for all $n\geq1$. Then,
\begin{equation}
0\geq v_{n}\left(  x_{n}\right)  +d_{\mathcal{O}^{c}}\left(  x_{n}\right)
\wedge1+\sup_{u\in U}\left\{
\begin{array}
[c]{c}%
-\left\langle f\left(  x_{n},u\right)  ,\nabla\varphi\left(  x_{n}\right)
+2\left(  x_{n}-x\right)  \right\rangle \\
-\lambda\left(  x_{n},u\right)  \int_{%
\mathbb{R}
^{N}}\left(  v_{n}(z)-v_{n}\left(  x_{n}\right)  \right)  Q\left(
x_{n},u,dz\right)
\end{array}
\right\}  .\label{A2.1'}%
\end{equation}
We have%
\begin{equation}
-\left\langle f\left(  x_{n},u\right)  ,\nabla\varphi\left(  x_{n}\right)
+2\left(  x_{n}-x\right)  \right\rangle \geq-\left\langle f\left(  x,u\right)
,\nabla\varphi\left(  x\right)  \right\rangle -C\left(  \left\vert
x_{n}-x\right\vert +\left\vert \nabla\varphi(x_{n})-\nabla\varphi
(x)\right\vert \right)  ,\label{A2.2}%
\end{equation}
where $C>0$ is a generic constant independent of $n\geq1$ and $u\in U$ which
may change from one line to another$.$ We also get
\begin{align}
- &  \lambda\left(  x_{n},u\right)  \int_{%
\mathbb{R}
^{N}}\left(  v_{n}(z)-v_{n}\left(  x_{n}\right)  \right)  Q\left(
x_{n},u,dz\right)  \nonumber\\
&  \geq-\lambda\left(  x,u\right)  \int_{%
\mathbb{R}
^{N}}\left(  v(z)-v\left(  x\right)  \right)  Q\left(  x,u,dz\right)
-C\left(  \left\vert x_{n}-x\right\vert +\left\vert v_{n}\left(  x_{n}\right)
-v\left(  x\right)  \right\vert +\eta_{v}\left(  \left\vert x_{n}-x\right\vert
\right)  \right)  \nonumber\\
&  -C\sup_{u\in U}\int_{%
\mathbb{R}
^{N}}\left\vert v_{n}\left(  z\right)  -v\left(  z\right)  \right\vert
Q\left(  x_{n},u,dz\right)  .\label{A2.3}%
\end{align}
Finally, for every $m\geq1,$
\begin{align}
&  \sup_{u}\int_{%
\mathbb{R}
^{N}}\left\vert v_{n}\left(  z\right)  -v\left(  z\right)  \right\vert
Q\left(  x_{n},u,dz\right)  \nonumber\\
&  \leq\sup_{z\in\overline{B}\left(  0,m+\left\vert x\right\vert +1\right)
}\left(  \left\vert v_{n}(z)-v(z)\right\vert \right)  +C\sup_{u\in U}Q\left(
x_{n},u,%
\mathbb{R}
^{N}\smallsetminus\overline{B}\left(  0,m+\left\vert x\right\vert +1\right)
\right)  \nonumber\\
&  \leq\sup_{z\in\overline{B}\left(  0,m+\left\vert x\right\vert +1\right)
}\left(  \left\vert v_{n}(z)-v(z)\right\vert \right)  +C\sup_{u\in U}Q\left(
x_{n},u,%
\mathbb{R}
^{N}\smallsetminus\overline{B}\left(  x_{n},m\right)  \right)  \nonumber\\
&  \leq\sup_{z\in\overline{B}\left(  0,m+\left\vert x\right\vert +1\right)
}\left(  \left\vert v_{n}(z)-v(z)\right\vert \right)  +C\sup_{y\in%
\mathbb{R}
^{N},u\in U}Q\left(  y,u,%
\mathbb{R}
^{N}\smallsetminus\overline{B}\left(  y,m\right)  \right)  .\label{A2.4}%
\end{align}
We substitute (\ref{A2.2})-(\ref{A2.4}) in (\ref{A2.1'}) and allow
$n\rightarrow\infty$ to have
\begin{align}
0 &  \geq v\left(  x\right)  +d_{\mathcal{O}^{c}}\left(  x\right)
\wedge1+\sup_{u\in U}\left\{  -\left\langle f\left(  x,u\right)
,\nabla\varphi\left(  x\right)  \right\rangle -\lambda\left(  x,u\right)
\int_{%
\mathbb{R}
^{N}}\left(  v(z)-v\left(  x\right)  \right)  Q\left(  x,u,dz\right)
\right\}  \nonumber\\
&  -C\sup_{y\in%
\mathbb{R}
^{N},u\in U}Q\left(  y,u,%
\mathbb{R}
^{N}\smallsetminus\overline{B}\left(  y,m\right)  \right)  ,
\end{align}
for all $m\geq1.$ We conclude using the Assumption A5.
\end{proof}


\begin{thebibliography}{99}                                                                                               %
\bibitem {AT}Alvarez, O., Tourin, A., \textit{Viscosity solutions of nonlinear
integro-differential equations,} Annales de l'institut Henri Poincar\'{e} (C)
Analyse non lin\'{e}aire (1996), 13 no. 3, pp. 293-317.

\bibitem {A}Aubin, J.-P., \textit{Viability Theory, }Birkh\"{a}user (1992).

\bibitem {ADP}Aubin, J.-P., Da Prato, G., \textit{Stochastic Viability and
invariance, }Annali Scuola Normale di Pisa (1990), No. 27, pp. 595-694.

\bibitem {AF}Aubin, J.-P., Frankowska, H., \textit{Set Valued Analysis},
Birkh\"{a}user, Boston (1990).

\bibitem {BC}M. Bardi and I. Capuzzo-Dolcetta, \textit{Optimal control and
viscosity solutions of Hamilton-Jacobi- Bellman equations}. Systems and
Control: Foundations and Applications, Birkh\"{a}user, Boston (1997).

\bibitem {BG}Bardi, M., Goatin, P., \textit{Invariant sets for controlled
degenerate diffusions: a viscosity solutions approach,} Stochastic analysis,
control, optimization and applications (1999), pp. 191--208, Systems Control
Found. Appl., Birkh\"{a}user Boston, Boston, MA.

\bibitem {BJ}Bardi, M., Jensen, R., \textit{A geometric characterization of
viable sets for controlled degenerate diffusions}, Set-Valued Anal. 10 (2002),
no. 2-3, pp. 129--141.

\bibitem {BI}Barles, G., Imbert, C., \textit{Second-Order Elliptic
Integro-Differential Equations: Viscosity Solutions Theory Revisited,} Annales
de l'IHP (2008), Vol. 25, No 3, pp. 567-585.

\bibitem {BaJ}Barles, G., Jakobsen, E.R., \textit{On the convergence rate of
approximation schemes for Hamilton-Jacobi-Bellman equations}, M2AN Math.
Model. Numer. Anal. 36 (2002), no. 1, pp. 33--54.

\bibitem {BGQ}Buckdahn, R., Goreac, D., Quincampoix, M.,\textit{ Stochastic
Optimal Control and Linear Programming Approach (}preprint, 2010).

\bibitem {BPQR}Buckdahn, R., Peng, S., Quincampoix, M., Rainer, C.,
\textit{Existence of stochastic control under state constraints, }C. R. Acad.
Sci. Paris S\'{e}r. I Math., 327 (1998), pp. 17-22.

\bibitem {CGT}Cook, D.,L., Gerber, A.,N., Tapscott, S.,J., \textit{Modelling
stochastic gene expression: Implications for haploinsufficiency, }Proc. Natl.
Acad. Sci. USA (1998), Vol. 95, pp. 15641-15646.

\bibitem {CDR}Crudu, A., Debussche, A., Radulescu, O., \textit{Hybrid
stochastic simplifications for multiscale gene networks}, BMC Systems Biology
(2009), 3: 89.

\bibitem {D}Davis, M.H.A.,\textit{ Markov Models and Optimization, }Monographs
on Statistics and Applied probability 49, Chapman \& Hall, London (1993).

\bibitem {De}Delbr\"{u}ck, M., \textit{Statistical Fluctuations in
Autocatalytic Reactions}, J. Chem. Phys. (1940), 8, pp. 120-124.

\bibitem {GT}Gautier, S., Thibault, L., \textit{Viability for constrained
stochastic differential equations, }Differential Integral Equations 6 (1993),
no. 6, pp. 1395-1414.

\bibitem {HPDC}Hasty, J., Pradines, J., Dolnik, M., Collins, J., J.,
\textit{Noise-based switches and amplifiers for gene expression, }PNAS (2000),
vol. 97, no. 5, pp. 2075-2080.

\bibitem {S}Soner, H.M., \textit{Optimal control with state-space constraint.
II.,} SIAM J. Control Optim. 24 (1986), no. 6, pp. 1110--1122. \textit{ }
\end{thebibliography}
\end{document}